\documentclass[11pt]{amsart}

\usepackage[foot]{amsaddr}
\usepackage{latexsym}
\usepackage{amsfonts}
\usepackage{amssymb}
\usepackage{amsmath}
\usepackage{mathrsfs}
\usepackage{bbm}
\usepackage{dsfont}
\usepackage{hyperref}
\allowdisplaybreaks

\newcommand{\be}{\begin{equation}}
\newcommand{\ee}{\end{equation}}
\newcommand{\bea}{\begin{eqnarray}}
\newcommand{\eea}{\end{eqnarray}}
\newcommand{\bean}{\begin{eqnarray*}}
\newcommand{\eean}{\end{eqnarray*}}
\newcommand{\brray}{\begin{array}}
\newcommand{\erray}{\end{array}}

\newtheorem{dfn}{Definition}[section]
\newtheorem{thm}[dfn]{Theorem}
\newtheorem{lmma}[dfn]{Lemma}
\newtheorem{ppsn}[dfn]{Proposition}
\newtheorem{crlre}[dfn]{Corollary}
\newtheorem{xmpl}[dfn]{Example}
\newtheorem{rmrk}[dfn]{Remark}

\newcommand{\bdfn}{\begin{dfn}\rm}
\newcommand{\bthm}{\begin{thm}}
\newcommand{\blmma}{\begin{lmma}}
\newcommand{\bppsn}{\begin{ppsn}}
\newcommand{\bcrlre}{\begin{crlre}}
\newcommand{\bxmpl}{\begin{xmpl}}
\newcommand{\brmrk}{\begin{rmrk}\rm}

\newcommand{\edfn}{\end{dfn}}
\newcommand{\ethm}{\end{thm}}
\newcommand{\elmma}{\end{lmma}}
\newcommand{\eppsn}{\end{ppsn}}
\newcommand{\ecrlre}{\end{crlre}}
\newcommand{\exmpl}{\end{xmpl}}
\newcommand{\ermrk}{\end{rmrk}}

\newcommand{\bbc}{\mathbb{C}}
\newcommand{\bbz}{\mathbb{Z}}

\newcommand{\bbn}{\mathbb{N}}
\newcommand{\bbr}{\mathbb{R}}

\newcommand{\bbt}{\mathbb{T}}

\author{ S. Sundar }
\title {On the  KMS states for the Bernoulli shift}

\begin{document}
\maketitle

\begin{abstract}
Let $\Omega:=\{0,1\}^{\bbz}$ be the Cantor space, and let $\tau:\Omega \to \Omega$ be the Bernoulli shift. 
 For the flow on the crossed product $C(\Omega)\rtimes_\tau \bbz$ determined by a potential that depends on only one coordinate, we show that for every $\beta \neq 0$, there is an extremal $\beta$-KMS state on $C(\Omega)\rtimes_\tau \bbz$ of type $II_\infty$. 
 Also, when the potential takes values that are rationally dependent,  we determine the values of $\lambda \in (0,1)$ for which there is a an extremal $\beta$-KMS state of type $III_\lambda$. 
\end{abstract}

 \noindent {\bf AMS Classification No. :} {Primary 46L30; Secondary 37A55.}  \\
{\textbf{Keywords :}} Bernoulli shift, KMS states, Isometric representations.

\section{Introduction}
The analysis of KMS states on various $C^{*}$-algebras has been in vogue for a long time and has  attracted quite a lot of attention  over the last two decades. We can consider \cite{Evans}, \cite{Cuntz_Laca_KMS}, \cite{Laca_Raeburn_KMS}, \cite{Laca_Neshveyev_KMS}, \cite{Thomsen_CMP}, \cite{Thomsen2021}, \cite{Thomsen}, \cite{Thomsen_Christensen1}, \cite{Kum_Ren} as a small sample of papers that deal with this subject.  One of the earliest important example
is the case of the Cuntz-algebra $\mathcal{O}_n$ which is the universal $C^{*}$-algebra generated by isometries $\{s_1,s_2,\cdots,s_n\}$ that satisfy the relation \[
\sum_{i=1}^{n}s_is_i^{*}=1.\]
Thanks to the universal property of $\mathcal{O}_n$, given positive real numbers $\lambda_1,\lambda_2,\cdots,\lambda_n$, there exists a $1$-parameter group of automorphisms $\sigma:=\{\sigma_t\}_{t \in \bbr}$ on $\mathcal{O}_n$ such that for every $k \in \{1,2,\cdots,n\}$ and $t \in \bbr$, $\sigma_t(s_k)=e^{it\lambda_k}s_k$.  Evans proved in \cite{Evans} that, with this action of $\bbr$, $\mathcal{O}_n$ has a unique KMS state attained at the inverse temperature $\beta$, where $\beta$ is the unique solution for the equation $\displaystyle \sum_{i=1}^{n}e^{-\beta \lambda_i}=1$. Izumi (\cite{Izumi_Cuntz}) showed that in the GNS representation, $\mathcal{O}_n$ generates a type III factor.

From the groupoid perspective, due to Renault, $\mathcal{O}_n$ is isomorphic to the $C^{*}$-algebra of the Deaconu-Renault groupoid associated to the one-sided  shift on the Cantor space $\{0,1,2,\cdots,n-1\}^{\bbn}$. In \cite{Kum_Ren}, Kumjian and Renault
showed that Evans' result can be deduced from this realisation, where the action of $\bbr$ is given by a cocycle which in turn is determined by a potential that depends only on the first coordinate. It is quite natural to ask what phenomena arise if we replace the one-sided shift by  the two-sided shift. 

 We consider only the two-sided shift on two symbols. The situation is vastly different.  As opposed to $\mathcal{O}_2$, the $C^{*}$-algebra $C(\{0,1\}^{\bbz})\rtimes \bbz$ is not simple. This could be one reason why the structure of  KMS states for the two-sided shift does not seem to be investigated in greater detail in the literature.  Neverthless, the author believes that there are a few questions that are of interest, and it 
is worth the effort to resolve them. We pose one such question and resolve it partially.

Let $\Omega:=\{0,1\}^{\bbz}$, and let $\tau$ be the Bernoulli shift on $\Omega$, i.e $\tau(x)_k=x_{k-1}$.  Let $\chi:\Omega \to \bbr$ be a continuous  function  that depends on only one coordinate.
The universal property of the crossed product $C(\Omega)\rtimes_{\tau}\bbz$ grants us a $1$-parameter group of automorphisms $\sigma:=\{\sigma_t\}_{t \in \bbr}$ on $C(\Omega)\rtimes \bbz$ such that for every $t \in \bbr$, 
$\sigma_t(f)=f$ for $f \in C(\Omega)$, and $\sigma_t(u)=ue^{it\chi}$. Here, $u$ is the canonical unitary of the crossed product $C(\Omega)\rtimes \bbz$. (This type of flow for a more general dynamical system $(X,T)$, where $X$ is a compact space and $T$ is a homeomorphism, was considered by Christensen and Thomsen in \cite{Thomsen}). 

Suppose that $\chi$ takes values $a$ and $b$. 
The analysis of KMS states is interesting only when $a$, $b$ are non-zero, and when they are of opposite signs. For a justification, the reader can consult Theorem 6.2 of \cite{Thomsen}, and the analysis carried out at   Page 18 in \cite{SAS}. Thus, we can normalise, and assume that $\chi: \Omega \to \bbr$ is given by
\begin{equation*}
\chi(x):=\begin{cases}
 1 & \mbox{ if
} x_{-1} =0,\cr
   &\cr
   - \theta &  \mbox{ if } x_{-1}=1
         \end{cases}
\end{equation*}
for some $\theta>0$. 
 The structure of KMS states for the flow $\sigma$, determined by the potential $\chi$, is vastly different from the $\mathcal{O}_2$ case. We have the following contrasting features. 
\begin{enumerate}
\item[(1)] The set of possible inverse temperatures is the whole real line $\bbr$.
\item[(2)] For every $\beta \in \bbr$ and for every $t \in \{I,II,III\}$, there is a continuum of extremal $\beta$-KMS states  of type $t$.
\end{enumerate}
In the case when $\theta$ is rational, it is not difficult to construct extremal KMS states of type $I$, $II$ and $III$. In the case when $\theta$ is irrational, it is still not hard to construct extremal KMS states
of type $I$ and type $II$. However, up to the author's knowledge, the construction of a type $III$  example (in fact, type $III_1$) requires the non-trivial work of Nakada (\cite{Nakada1}, \cite{Nakada}). 
The type $II$ examples constructed in \cite{SAS} are of type $II_1$. For more details, we refer the reader to \cite{SAS} (see also Remark \ref{continuum of states}).

Given that all the three factorial KMS states are possible, a more refined question would be to  ask whether  every possible Krieger type occurs or not. In particular, the author believes that it is worth  asking the following  questions.
Let $\beta$ be a non-zero real number. 
\begin{enumerate}
\item[(1)] Does there exist an extremal $\beta$-KMS state of type $II_\infty$?
\item[(2)] Determine the values of $\lambda \in (0,1)$ for which there is an extremal $\beta$-KMS state for $\sigma$ of type $III_\lambda$.
\item[(3)] Does there exist an extremal $\beta$-KMS state of type $III_0$? 
\end{enumerate}
As mentioned earlier, it follows from the work of Nakada (\cite{Nakada1}, \cite{Nakada}) that type $III_1$ occurs when $\theta$ is irrational.  In this paper, we show the existence of a type $II_\infty$ example. 
Also, when $\theta$ is rational, we determine the values of $\lambda \in (0,1)$ for which there is a $\beta$-KMS state of type $III_\lambda$. 

Elementary considerations establish   a 
bijective correspondence 
$
m \to \omega_m$
between the set of  non-atomic, ergodic probability measures $m$ on $\Omega$ that are $e^{-\beta \chi}$-conformal, i.e. \[\frac{d(m\circ \tau)}{dm}=e^{-\beta \chi}\]
and the set of extremal $\beta$-KMS states on $C(\Omega)\rtimes_\tau \bbz$ that are not of type $I$. Moreover, for every $t \in \{II_1,II_\infty, III_\lambda: \lambda \in [0,1]\}$, $m$ is of type $t$ if and only if $\omega_m$ is of type $t$.
Our  main theorem is stated below.

\begin{thm}
\label{main}
 For every non-zero $\beta \in \bbr$, there exists an ergodic  probability measure $m$ on $\Omega$ such that $m$ is $e^{-\beta \chi}$-conformal, and $m$ is of type $II_\infty$. 
\end{thm}

The proof of the above  theorem is existential and not constructive. Our proof is also  operator algebraic. Our first reduction, which we undertake in Section 2,  is based on the bijection established in \cite{SAS} that asserts that the  set of $e^{-\beta \chi}$-conformal measures on $\Omega$ is in 
bijective correspondence with the set of $e^{-\beta c}$-conformal, Radon measures on the unit space of the Deaconu-Renault  groupoid $X_u \rtimes \bbn^2$  that encodes two parameter discrete semigroups of isometries with commuting range projections. Here, the cocycle $c:X_u \rtimes \bbn^2 \to \bbr$ is given by the homomorphism $c:\bbz^2 \to \bbr$ defined by $c(m,n)=m+n\theta$.

Secondly, we show that constructing $e^{-\beta c}$-conformal measures on $X_u$ is equivalent to constructing representations ( equivalently, constructing  semigroups  of isometries, indexed by $\bbn^2$, with commuting range projections) of $C^{*}(X_u\rtimes \bbn^2)$ for which the eigenspace of $\bbn^2$ corresponding to the character $e^{-\frac{\beta c}{2}}$ is non-trivial. This translation is inspired by the techniques of \cite{Wu}. We undertake this translation in Section 3, where we prove the results in a more general setting of a closed subsemigroup of a locally compact, abelian group. We show that proving Thm. \ref{main} amounts to  constructing a two parameter semigroup of isometries with commuting range projections satisfying certain hypotheses. In Section 4, by appealing to certain results available in the literature concerning ergodic theory, we produce  such a semigroup of isometries which proves Thm. \ref{main}.

As far as  type $III$ examples are concerned, we prove the following. 
Suppose $\theta$ is rational, and suppose $\theta=\frac{p}{q}$ with $\gcd(p,q)=1$. Let $\beta$ be a non-zero real number. Suppose $m$ is an $e^{-\beta \chi}$-conformal measure that is of type $III_\lambda$. 
 Since $\chi$ takes values in $\frac{1}{q}\bbz$, it follows that the ratio set \[r(\tau)\cap (0,\infty) \subset  \{e^{-\frac{\beta n}{q}}: n \in \bbz\}\] which is a closed subgroup of $(0,\infty)$. Hence, $\lambda$ is necessarily of the form $e^{-\frac{|\beta| n}{q}}$ for some $n \geq 1$. 
The next theorem ensures that every such $\lambda$ is realised. 

\begin{thm}
\label{main1}
Suppose $\theta:=\frac{p}{q}$ is rational, and let $\beta$ be a non-zero real number. Then, for every $\lambda \in \{e^{-\frac{|\beta| n}{q}}: n \in \{1,2,\cdots\}\}$, there exists an ergodic, probability measure $m$ on $\Omega$ such that $m$ is $e^{-\beta \chi}$-conformal, and $m$ is  of type $III_\lambda$.
\end{thm}

The author's interest to analyse the structure of KMS states on $C(\Omega)\rtimes_\tau \bbz$ stems from the close relationship, established in \cite{SAS}, that exists between the KMS states on $C_{c}^{*}(\bbn^2)$ (which is the universal $C^{*}$-algebra generated by a semigroup of isometries $\{v_{(m,n)}\}_{(m,n) \in \bbn^2}$ with commuting range projections) for the time evolution determined by the homomorphism $c:\bbz^2 \to \bbr$ and the KMS states on $C(\Omega)\rtimes \bbz$.  Roughly, at positive inverse temperature, every KMS state on $C_{c}^{*}(\bbn^2)$ is an 'amplified version' of a unique KMS state on $C(\Omega) \rtimes \bbz$. Conversely, every KMS state on $C(\Omega)\rtimes \bbz$ is a `corner' of a unique KMS state on $C_c^{*}(\bbn^2)$. From this, it is apparent that type $II$ KMS states on $C_{c}^{*}(\bbn^2)$ are always of type $II_\infty$. But it is not at all clear why type $II_\infty$ occurs for $C(\Omega)\rtimes \bbz$. In fact, the examples constructed in \cite{SAS} are of type $II_1$. Thus, it 
is of intrinsic interest to investigate whether there are type $II_\infty$ KMS states on $C(\Omega)\rtimes \bbz$. 

\textit{Convention:} For us, the set of natural numbers $\bbn$ contains $0$. All the Hilbert spaces considered in this paper are assumed to be separable, and the inner product is linear in the first variable.

\section{A reduction}
Let $\Omega:=\{0,1\}^{\bbz}$ be the Cantor space, and let $\tau:\Omega \to \Omega$ be the Bernoulli shift defined by $\tau(x)_k=x_{k-1}$. Suppose $\theta>0$. Let $\chi:\Omega \to \bbr$ be the continuous map defined by 
\begin{equation*}
\chi(x):=\begin{cases}
 1 & \mbox{ if
} x_{-1} =0,\cr
   &\cr
   - \theta &  \mbox{ if } x_{-1}=1.
         \end{cases}
\end{equation*}
Recall that the potential $\chi$ defines a flow $\sigma:=\{\sigma_t\}_{t \in \bbr}$ on $C(\Omega)\rtimes_\tau \bbz$, where the automorphism $\sigma_t$ is given by 
\[
\sigma_t(f)=f~~\textrm{and~~}\sigma_t(u)=ue^{it\chi}\]
for $f \in C(\Omega)$. Here, $u$ stands for the canonical unitary of the crossed product $C(\Omega)\rtimes_\tau \bbz$. 

 Fix a real number $\beta$. Let $\omega$ be a $\beta$-KMS state on $C(\Omega)\rtimes_\tau \bbz$ for $\sigma$. Let $m:=m_\omega$ be the probability measure on $\Omega$ that corresponds to the state $\omega|_{C(\Omega)}$. Then, $m$ is an $e^{-\beta \chi}$-conformal measure on $\Omega$, i.e. 
\[
m(\tau(B))=\int_{B}e^{-\beta \chi}dm\]
for every Borel set $B \subset \Omega$. 
Conversely, suppose $m$ is an $e^{-\beta \chi}$-conformal probability measure on $\Omega$. Then, there exists a unique  $\beta$-KMS state $\omega=\omega_m$ on $C(\Omega)\rtimes_\tau \bbz$ such that 
\[\omega(fu^k)=\delta_{k,0}\int f dm.\] 
Also, $m_{\omega_m}=m$. 

For a $\beta$-KMS state $\omega$ on $C(\Omega)\times_\tau \bbz$, we denote the associated GNS representation by $\pi_\omega$. Let us make the following observations regarding the set of extremal $\beta$-KMS states. 
\begin{enumerate}
\item[(1)] Suppose $m$ is a non-atomic probability measure on $\Omega$. Since $\Omega$ has only countably many periodic points, it follows that the action of $\bbz$ on $(\Omega,m)$ via $\tau$ is essentially free. 
\item[(2)] Suppose $m$ is an $e^{-\beta \chi}$-conformal probability measure that is non-atomic and ergodic. It follows from  Thm. 4.13 of \cite{Thomsen} that $\omega_m$ is an extremal $\beta$-KMS state
               of type $t$ where $t \in \{II, III\}$. 
 \item[(3)] Suppose $\omega$ is an extremal $\beta$-KMS state. Let $m$ be the associated $e^{-\beta \chi}$-conformal measure. Then, it follows from Lemma 3.7 of \cite{Thomsen} that $m$ is ergodic. 
 \begin{enumerate}
 \item[(a)]
 Suppose $m$ is atomic and concentrated on an orbit of a non-periodic point. Then, by Corollary 1.4 of \cite{Neshveyev}, it follows that $\omega=\omega_m$. In this case, it follows from Thm. 4.13 of \cite{Thomsen} that $\omega$ is of type $I_\infty$. 
 \item[(b)]  Suppose $m$ is atomic and concentrated on an orbit of a periodic point of period $p$. Then, by Corollary 1.4 of \cite{Neshveyev}, it follows that there exists $z \in \bbt$ such that 
  \begin{equation*}
\omega(f u^k):=\begin{cases}
 0  & \mbox{ if
} k \notin p\bbz \cr
   &\cr
    z^k\displaystyle \int_{\Omega} f(x)dm(x) &  \mbox{ if $k \in p\bbz$}.
         \end{cases}
\end{equation*}
In this case, it is not difficult to prove  as in Prop. 3.4 of \cite{SAS} that the von Neumann algebra $\pi_\omega(C(\Omega)\rtimes_\tau \bbz)^{''}=L^{\infty}(\Omega,m)\rtimes \bbz/p\bbz \cong M_p(\bbc)$. In this case, $\omega$ is of type $I_p$.
\item[(c)] Suppose $m$ is non-atomic. Since $\Omega$ has only countably many periodic points and $m$ is non-atomic, it follows from  Thm. 1.3 of \cite{Neshveyev} that $\omega=\omega_m$. In this case, it follows from Thm. 4.13 of \cite{Thomsen} that $\omega$ is of type $t$ where $t \in \{II,III\}$. 
\end{enumerate}
\end{enumerate}

Thus, the association
\[
m \to \omega_m\]
sets up a bijective correspondence between  the set of non-atomic, ergodic, $e^{-\beta \chi}$-conformal  probability measures on $\Omega$ and the set of extremal $\beta$-KMS states on $C(\Omega)\rtimes_\tau \bbz$  that are not of type $I$. 
Thanks to Thm. 4.13 of \cite{Thomsen}, for an ergodic, $e^{-\beta \chi}$-conformal probability measure $m$ on $\Omega$ and $t \in \{II_1,II_\infty,III_\lambda\}$, $m$ is of type $t$ if and only if $\omega_m$ is of type $t$.
Thus, exhibiting an extremal $\beta$-KMS state of type $II_\infty$ is equivalent to exhibiting an ergodic, $e^{-\beta \chi}$-conformal probability measure $m$ on $\Omega$ such that $m$ is of type $II_\infty$. 

For a real number $\beta$, 
let $\mathcal{M}_{\beta}(\Omega)$ denote the set of $e^{-\beta \chi}$-conformal probability measures on $\Omega$, and let 
\[
\mathcal{M}_{e,\beta}(\Omega):=\{ m \in \mathcal{M}_{\beta}(\Omega): \textrm{$m$ is ergodic}\}.\]

Next, we establish a bijection between $\mathcal{M}_{\beta}(\Omega)$ and $\mathcal{M}_{-\beta}(\Omega)$ for every $\beta$. Let $\kappa:\Omega \to \Omega$ be the homeomorphism defined by
\[
\kappa(x)_k=x_{-k-1}.\]
Note that $\kappa$ is of order two, and $\tau \circ \kappa=\kappa \circ \tau^{-1}$. 

\begin{ppsn}
\label{symmetry}
Let $\beta \in \bbr$ be given. If $m \in \mathcal{M}_{\beta}(\Omega)$, then $m \circ \kappa \in \mathcal{M}_{-\beta}(\Omega)$. Moreover, the map
\[
\mathcal{M}_{\beta}(\Omega) \ni m \to m \circ \kappa \in \mathcal{M}_{-\beta}(\Omega)\]
is a bijection. Also, the map 
\[
\mathcal{M}_{\beta}(\Omega) \ni m \to m \circ \kappa \in \mathcal{M}_{-\beta}(\Omega)\]
preserves ergodicity and the Krieger type. 
\end{ppsn}
\textit{Proof.}  Using the fact that $\tau \circ \kappa=\kappa \circ \tau^{-1}$, it is routine to check that for a probability measure $m$ on $\Omega$, $m \in \mathcal{M}_{\beta}(\Omega)$ if and only if $m \circ \kappa \in \mathcal{M}_{-\beta}(\Omega)$. 
Again, for $m \in \mathcal{M}_{\beta}(\Omega)$, using the fact that $\tau \circ \kappa=\kappa \circ \tau^{-1}$, it can be verified that $m$ is ergodic if and only if $m \circ \kappa$ is ergodic.

Let $m \in \mathcal{M}_{e,\beta}(\Omega)$ be given. Clearly, $m$ is atomic if and only if $m \circ \kappa$ is atomic. 
Suppose that $m$ is non-atomic. Idenfity $L^{\infty}(\Omega, m \circ \kappa)$ with $L^{\infty}(\Omega,m)$ via the map 
\[
L^{\infty}(\Omega,m \circ \kappa) \ni f \to f \circ \kappa \in L^{\infty}(\Omega,m).\]
With the above identification, and thanks to the equality $\kappa \circ \tau=\tau^{-1} \circ \kappa$, we have 
\[
L^{\infty}(\Omega, m \circ \kappa) \rtimes_{\tau} \bbz \cong L^{\infty}(\Omega,m) \rtimes_{\tau^{-1}}\bbz \cong L^{\infty}(\Omega,m) \rtimes_{\tau} \bbz.\]
Therefore, $m$ and $m \circ \kappa$ have the same Krieger type. This completes the proof. \hfill $\Box$

The first key step is to convert the problem of constructing conformal measures on $\Omega$ into a problem of constructing conformal measures on  `the universal dynamical system' that encodes discrete two parameter
semigroups of isometries with commuting range projections. First, we make a few definitions in a more general context. 

Let $G$ be a locally compact, second countable, Hausdorff abelian group. A second countable, locally compact, Hausdorff space which has a continuous $G$-action will be called a $G$-space. Let $Y$ be a  $G$-space, and suppose $\mu$ is a quasi-invariant, ergodic, non-zero, Radon measure on $Y$. We say that 
\begin{enumerate}
\item[(1)] $\mu$ is of type $I$ if $\mu$ is supported on an orbit.
\item[(2)] $\mu$ is of type $II_1$ if $\mu$ is not of type $I$, and if the measure class $[\mu]$ has a  a non-zero $G$-invariant, Radon measure\footnote{We warn the reader that,  type $II_1$, in this sense, does not mean that the associated crossed product $L^{\infty}(Y,\mu)\rtimes G$ is a type $II_1$ factor.}.
\item[(3)] $\mu$ is of type $II_\infty$ if $\mu$ is not of type $I$, and if the measure class $[\mu]$ has a non-zero $\sigma$-finite, $G$-invariant measure but has no non-zero $G$-invariant, Radon measure.
\item[(4)] $\mu$ is of type $III$ if the measure class $[\mu]$ has no non-zero $\sigma$-finite, $G$-invariant measure. 
\end{enumerate}
If $\mu$ is of type $III$ and $\lambda \in [0,1]$, we say that $\mu$ is of type $III_\lambda$ if $L^{\infty}(Y_u)\rtimes G$ is a factor of type $III_\lambda$. 

Suppose $c:G \to \bbr$ is a continuous homomorphism, and suppose $Y$ is a $G$-space. A non-zero, Radon measure  $\mu$ is said to be $e^{-\beta c}$-conformal if $\mu(E+s)=e^{-\beta c(s)}\mu(E)$ for every $s \in G$ and for every Borel subset $E \subset Y$. 

Let $G$ be a countable, discrete abelian group, and let $P \subset G$ be a subsemigroup containing $0$ such that $P-P=G$. Define 
\begin{align*}
\overline{Y}_u&:=\{A \subset G: -P+A \subset A, A \neq \emptyset\},\\
\overline{X}_u&:=\{A \in \overline{Y}_u: 0 \in A\}.
\end{align*}
We identify $\overline{Y}_u$ with a subset of $\{0,1\}^{G}$ in the usual way, and we endow $\overline{Y}_u$ with the subspace topology inherited from the product topology
on $\{0,1\}^{G}$. Note that $\overline{Y}_u$ is a locally compact, Hausdorff space, and $\overline{X}_u$ is a compact subset of $\overline{Y}_u$. The map 
\[
\overline{Y}_u \times G \ni (A,s) \to A+s \in \overline{Y}_u\]
defines an action of $G$ on $\overline{Y}_u$, and $P$ leaves $\overline{X}_u$ invariant, i.e. $\overline{X}_u +P \subset \overline{X}_u$. 

\begin{rmrk}
It was demonstrated in \cite{Sundar_Ore} that the dynamical system $(\overline{Y}_u,G)$  encodes all semigroups of isometries, indexed by $P$, that has commuting range projections. In particular, $C^{*}(\overline{X}_u \rtimes P)$, where $\overline{X}_u \rtimes P$ is the Deaconu-Renault groupoid, is the universal $C^{*}$-algebra generated by isometries $\{v_a:a \in P\}$ such that 
\begin{enumerate}
\item[(1)] for $a,b \in P$, $v_{a+b}=v_av_b$, and
\item[(2)] the family $\{v_av_{a}^{*}:a \in P\}$ is a commuting family of projections.
\end{enumerate}
Moreover, $C^{*}(\overline{X}_u \rtimes P)$ is a full corner in $C_0(\overline{Y}_u) \rtimes G$. For more details, the reader is referred to \cite{Sundar_Ore} and Section 2 of \cite{SAS}.
\end{rmrk}

Set $Y_u:=\overline{Y}_u\backslash \{G\}$, and let $X_u:=\overline{X}_u \cap Y_u$. In \cite{SAS}, when $G=\bbz^2$ and $P=\bbn^2$, a nice parametrisation of $Y_u$ was obtained which we explain next. 
For the remainder of this section, assume $G=\bbz^2$ and $P=\bbn^2$. Set $e_1:=(1,0)$, $e_2:=(0,1)$, $v_1:=e_1$, and $v_2:=e_1+e_2$. Define a $\bbz^2$-action on $\Omega \times \bbz$ by
setting
\begin{align*}
(x,t)+v_1:&=(\tau(x),x_{-1}+t),\\
(x,t)+v_2:&=(x,t+1).
\end{align*}
For $(x,t) \in \Omega \times \bbz$, let $a(x,t)$ be the bi-infinite sequence  defined by 
\begin{equation}
 a(x,t)_m:=\begin{cases}
  t-(x_0+x_1+\cdots +x_{m-1}) & \mbox{if~} m> 0, \cr
      t & \mbox{if $m=0$~}, \cr
    t+(x_{-1}+x_{-2}+\cdots +x_{m}) & \mbox{if~} m<0.
         \end{cases}
\end{equation}
For $(x,t) \in \Omega \times \bbz$, set 
\[
A(x,t):=\{mv_1+nv_2: n \leq a(x,t)_m\}.\]
Then, $A(x,t) \in Y_u$ for every $(x,t) \in \Omega \times \bbz$. The following proposition is Prop. 4.1 of \cite{SAS}.
\begin{ppsn}[\cite{SAS}]
\label{parametrisation of $Y_u$}
The map 
\[
\Omega \times \bbz \ni (x,t) \to A(x,t) \in Y_u\]
is a $\bbz^2$-equivariant homeomorphism. 
\end{ppsn}
We identify $Y_u$ with $\Omega \times \bbz$ via the map given by Prop. \ref{parametrisation of $Y_u$}, and by abusing notation, we write $Y_u=\Omega \times \bbz$. Then, $X_u=\Omega \times \bbn$. 
Let $c:\bbz^2 \to \bbr$ be the homomorphism such that 
$c(e_1)=1$ and $c(e_2)=\theta$.  Let $\beta \in \bbr$. Denote the set of $e^{-\beta \chi}$-conformal probability measures on $\Omega$ by $\mathcal{M}_{\beta}(\Omega)$, and denote the set of $e^{-\beta c}$-conformal, non-zero, Radon measures on $Y_u$
by $\mathcal{M}_{\beta}(Y_u)$. Set 
\begin{align*}
\mathcal{M}_{e,\beta}(\Omega):&=\{m \in \mathcal{M}_{\beta}(\Omega): \textrm{$m$ is ergodic}\},\\
\mathcal{M}_{e,\beta}(Y_u):&=\{\mu \in \mathcal{M}_\beta(Y_u): \textrm{$\mu$ is ergodic}\},\\
\mathcal{M}_{\beta,v_2}(Y_u):&=\{\mu \in \mathcal{M}_\beta(Y_u): \mu(X_u\backslash(X_u+v_2))=1\},\\
\mathcal{M}_{e,\beta,v_2}(Y_u):&=\mathcal{M}_{e,\beta}(Y_u) \cap \mathcal{M}_{\beta,v_2}(Y_u).
\end{align*}

For a probability measure $m$ on $\Omega$, define a Radon measure $\overline{m}$ on $Y_u=\Omega \times \bbz$ by setting
\begin{equation}
\label{association}
\overline{m}(E\times \{n\})=e^{-n\beta(1+\theta)}m(E)\end{equation}
where $E \subset \Omega$ is a Borel subset. 

It is routine to prove the following proposition as in Prop. 4.2 of \cite{SAS}. Hence, we omit the proof. 
\begin{ppsn}
\label{bijection between conformal measures}
The map 
\[
\mathcal{M}_{\beta}(\Omega) \ni m \to \overline{m} \in \mathcal{M}_{\beta,v_2}(Y_u)\]
is a bijection. Moreover, for $m \in \mathcal{M}_{\beta}(\Omega)$, $m$ is ergodic if and only if $\overline{m}$ is ergodic. Suppose $m \in \mathcal{M}_{e,\beta}(\Omega)$.  Then, for $t \in \{I,II_1,II_\infty,III\}$, $m$ is of type $t$ if and only if $\overline{m}$ is of type $t$.
\end{ppsn} 

\begin{rmrk}
\label{continuum of states}

Let $\mu \in \mathcal{M}_\beta(Y_u)$. It is clear from the condition $$\mu(\Omega+\{n+1\})=\mu((\Omega \times \{n\})+v_2))=e^{-\beta(1+\theta)}\mu(\Omega \times \{n\})$$ that 
$\mu(\Omega \times \{0\}) \neq 0$. Thus, there exists a unique positive real number $r$ and a unique $m \in \mathcal{M}_{\beta}(\Omega)$ such that $\mu=r\overline{m}$.

For $\beta>0$, it was proved in \cite{SAS} that there is a continuum of  $e^{-\beta c}$-conformal, ergodic, non-zero, Radon measures on $Y_u$ of type $t$ for each $t \in \{I,II,III\}$. Combining Prop. \ref{bijection between conformal measures} and Prop. \ref{symmetry}, we see that for every non-zero $\beta$, there are uncountably many ergodic, $e^{-\beta \chi}$-conformal probability measures on $\Omega$ of type $t$ for each $t \in \{I,II,III\}$. 

\end{rmrk}

\begin{lmma}
\label{krieger type}
Let $m \in \mathcal{M}_{\beta}(\Omega)$ be such that $m$ is ergodic and is of type $III$. Let $\overline{m}$ be the Radon measure on $Y_u$ defined by Eq. \ref{association}.  For $\lambda \in [0,1]$, $m$ is of type $III_\lambda$ if and only if $\overline{m}$ is of type $III_\lambda$.  
\end{lmma}
\textit{Proof.} Note that $\overline{m}$ on $Y_u=\Omega \times \bbz$ is absolutely continuous w.r.t. $dmdn$, where $dn$ is the counting measure on $\bbz$. 
The crossed product $L^{\infty}(Y_u) \rtimes \mathbb{Z}^{2}$ can be written as an iterated crossed product $(L^{\infty}(Y_u) \rtimes \mathbb{Z}v_2) \rtimes \mathbb{Z}v_1$. But \[L^{\infty}(Y_u) \cong L^{\infty}(\Omega,m) \otimes \ell^{\infty}(\mathbb{Z}).\] Note that  the $v_2$-action is on the second factor, and it acts by translation by $1$ on $\mathbb{Z}$. Therefore, $L^{\infty}(Y_u) \rtimes \mathbb{Z}v_2 \cong L^{\infty}(\Omega ,m) \otimes B(\ell^{2}(\mathbb{Z}))$. 

Let $\tau:L^{\infty}(\Omega,m) \to L^{\infty}(\Omega,m)$ be the map defined by  $\tau(f)(x)=f(\tau^{-1}x)$.
   Denote the bilateral shift on $\ell^{2}(\mathbb{Z})$ by $U$, i.e. $Ue_n=e_{n+1}$, where $\{e_n: n \in \mathbb{Z}\}$ is the standard orthonormal basis for $\ell^2(\bbz)$. Define $p \in L^{\infty}(\Omega,m)$ by $p:=1_{\{x_0=0\}}$.

Once $L^{\infty}(Y_u) \rtimes \mathbb{Z}v_2$ is identified with $L^{\infty}(\Omega,m) \otimes B(\ell^{2}(\mathbb{Z}))$, by a routine direct computation, we see that the $v_1$-action on $L^{\infty}(\Omega.m) \otimes B(\ell^{2}(\mathbb{Z}))$ coincides with the automorphism 
$Ad(p \otimes 1+ (1-p) \otimes U) \circ (\tau \otimes 1)$.  Thus, the automorphism corresponding to the $v_1$-action is outer conjugate to $\tau \otimes 1$. Hence,
\[
L^{\infty}(Y_u,\overline{m})\rtimes \bbz^2 \cong (L^{\infty}(\Omega,m)\rtimes \bbz)\otimes B(\ell^2(\bbz)).\]
The result follows. \hfill $\Box$

We have now proved that proving Thm. \ref{main}  is equivalent to exhibiting an  ergodic, non-zero, Radon measure $\mu$ on $Y_u$ such that $\mu$ is $e^{-\beta c}$-conformal, and $\mu$ is of type $II_\infty$. How to construct ergodic, $e^{-\beta c}$-conformal, Radon measures on $Y_u$? 
We take inspiration from the results of \cite{Wu}. Let us recall the 
main result of \cite{Wu}. Suppose $\Gamma$ is a discrete group that acts on a compact space $X$ by homeomorphisms. Then, the set of $\Gamma$-invariant,
ergodic probability measures on $X$ is in bijective correspondence with the collection (up to unitary equivalence) of irreducible representations of the 
crossed product $C(X)\rtimes \Gamma$ for which $\Gamma$ has an invariant unit vector. We establish an analogous result here. We show that 
$\mathcal{M}_{e,\beta,v_2}(Y_u)$ is in bijective correspondence with the set of irreducible representations of $C^{*}(\overline{X}_u \rtimes \bbn^2)$ for 
which the eigenspace of $\bbn^2$ corresponding to the character $e^{-\frac{\beta c}{2}}$ is non-trivial. In the next section, we prove this bijection for a general  closed, subsemigroup of a locally compact  abelian group.

\begin{rmrk}
The author's reason for proving the results in the topological setting, and not just in the discrete setting, is because the results developed in the next section are needed for future applications to $E_0$-semigroups.
To avoid duplication, in the next section, the results are developed in the setting of subsemigroups of locally compact abelian groups.
\end{rmrk}


\section{A translation}
Let $G$ be a locally compact, second countable, Hausdorff abelian group. Let $P$ be a closed subsemigroup of $G$ containing $0$. We assume that $P-P=G$ and $P$ has dense interior. For $s,t \in G$, we say $s \leq t$ if $t-s \in P$, and $s<t$ if $t-s \in Int(P)$. 
We assume that $P$ has an order unit, i.e. there exists $a_0 \in Int(P)$ such that for every $s \in G$, there exists a positive integer $n$ such that $na_0>s$. Let $c:G \to \bbr$ be a continuous group homomorphism. The group $G$, the semigroup $P$, and the homomorphism $c:G \to \bbr$ are fixed for the rest of this section. 

We first collect a few definitions concerning semigroups of isometries. 
A semigroup of isometries $V:=\{V_a\}_{a \in P}$ is said to be \emph{pure} if $\bigcap_{a \in P}Ran(V_a)=\{0\}$. We also call a strongly continuous, semigroup of isometries, indexed by $P$, an isometric representation of $P$. Let $V:=\{V_a\}_{a \in P}$ be an isometric representation of $P$ on a Hilbert space $H$. Let $U:=\{U_s\}_{s \in G}$ be a strongly continuous, group of unitary operators on a Hilbert space $K$. We say that $(U,K)$ is the \emph{minimal unitary dilation} of $(V,H)$ if 
\begin{enumerate}
\item[(1)] the Hilbert space $K$ contains $H$ as a closed subspace, 
\item[(2)] for $a \in P$, $U_a|_{H}=V_a$, and
\item[(3)] the union $\bigcup_{a \in P}U_{a}^{*}H$ is dense in $K$.
\end{enumerate}
For the existence and the uniqueness (up to a unitary equivalence) of the minimal unitary dilation of an isometric representation, we refer the reader to \cite{Laca95}.

Let $V:=\{V_{a}\}_{a \in P}$ be a pure, isometric representation of $P$ on a Hilbert space $H$. Set 
\[
\mathcal{D}_{V}:=\bigcup_{a \in P}Ker(V_a^*).\]
Since $V$ is pure, $\mathcal{D}_V$ is a dense subspace of $H$. Also, for every $a \in P$, $V_a$ leaves $\mathcal{D}_V$ invariant. Let 
\[
\mathcal{L}_{V}:=\{\phi: \mathcal{D}_V \to \bbc: \textrm{$\phi$ is linear, and $\phi|_{Ker(V_a^*)}$ is bounded for every $a \in P$}\}.\]
For $a \in P$, let $T_a:\mathcal{L}_V \to \mathcal{L}_V$  be defined by 
\[
T_a\phi(\xi):=\phi(V_a\xi).\]
Then, $T:=\{T_a\}_{a \in P}$ is a semigroup of linear operators on $\mathcal{L}_V$. 

\begin{dfn}
Let $\beta \in \bbr$. With the foregoing notation, we say that $V$ is $e^{-\frac{\beta c}{2}}$-conformal if there exists a non-zero element $\phi \in \mathcal{L}_V$ such that 
\[
T_a\phi=e^{-\frac{\beta c(a)}{2}}\phi\]
for every $a \in P$. 
We say that $V$ is $1$-conformal if $V$ is $e^{-\frac{\beta c}{2}}$-conformal for $\beta=0$. For $\beta \in \bbr$, let \[
\mathcal{A}_\beta(V):=\{\phi \in \mathcal{L}_V: T_a\phi=e^{-\frac{\beta c(a)}{2}}\phi ~~\textrm{for all $a \in P$}\}.
\]

\end{dfn}

\begin{rmrk}
\label{additivity of index}
Suppose $V=W\oplus W^{'}$. Then, for $\beta \in \bbr$, we have \[\mathcal{A}_{\beta}(V)=\mathcal{A}_{\beta}(W)\oplus \mathcal{A}_{\beta}(W^{'}).\] Thus,
\[
\dim \mathcal{A}_\beta(V)=\dim \mathcal{A}_\beta(W) + \dim \mathcal{A}_\beta (W^{'}).\] 
Hence, if $V$ is not $e^{-\frac{\beta c}{2}}$-conformal, then $W$ is not $e^{-\frac{\beta c}{2}}$-conformal. 

\end{rmrk}

\begin{lmma}
\label{eigen vector implies conformality}
Keep the foregoing notation. Let $\beta \in \bbr$. Suppose there exists a non-zero vector $\xi \in H$ such that \[V_a^{*}\xi=e^{-\frac{\beta c(a)}{2}}\xi\] for all $a \in P$.
Then, $V$ is $e^{-\frac{\beta c}{2}}$-conformal. 
\end{lmma}
\textit{Proof.} Define $\phi:\mathcal{D}_V \to \bbc$ by $\phi(\eta):=\langle \eta|\xi\rangle$. Since $\mathcal{D}_V$ is dense in $H$ and $\xi \neq 0$, $\phi \neq 0$. Clearly,  $\phi \in \mathcal{L}_V$, and 
\[T_a\phi=e^{-\frac{\beta c(a)}{2}}\phi\]for every $a \in P$. The proof is complete. \hfill $\Box$

\begin{rmrk}
\label{coherence}
Let $V:=\{V_a\}_{a \in P}$ be a pure, strongly continuous semigroup of isometries on a Hilbert space $H$. For $a \in P$, let $E_a:=V_{a}V_{a}^{*}$. 
 Let $\phi \in \mathcal{L}_V$. For $a \in P$, let $\xi_a \in Ker(V_a^*)$ be such that $\phi(\xi)=\langle \xi|\xi_a\rangle$ for $\xi \in Ker(V_a^*)$. Then, the section \[P \ni a \to \xi_a \in \coprod_{b \in P}Ker(V_b^*)\] is coherent, i.e. 
given $a,b \in P$ with $a \leq b$, we have $E_{a}^{\perp}\xi_b=\xi_a$. 

Conversely, given a coherent section $\displaystyle \xi:P \to \coprod_{a \in P}Ker(V_a^*)$, there exists a unique $\phi \in \mathcal{L}_V$ such that, given $a \in P$, $\phi(\xi)=\langle \xi|\xi_a \rangle$  for every $\xi \in Ker(V_a^*)$. 
This way, we can identify $\mathcal{L}_V$ with the set of coherent sections $\displaystyle \xi:P \to \coprod_{a \in P}Ker(V_a^*)$. We freely use this identification. 

The space $\mathcal{L}_V$ can be given a Fr\'echet space structure, and the space $\mathcal{L}_V$ appears naturally in the analysis of decomposable product systems (\cite{Arveson}, \cite{Namitha_Sundar}) in the theory of $E_0$-semigroups.
\end{rmrk}

 Let $V:=\{V_a\}_{a \in P}$ be an isometric representation of $P$ on a Hilbert space $H$. For $a \in P$, set $E_a:=V_aV_a^{*}$. We say that $V$ has commuting range projections 
if $\{E_a:a \in P\}$ is a commuting family of projections. We say that $V$ is \emph{irreducible} if $\{V_a,V_a^{*}:a \in P\}^{'}=\bbc$.

Next, we recall from \cite{Sundar_Ore} the relevant results that we need for this paper. Fix an order unit $a_0 \in Int(P)$ once and for all for the rest of this section. 
Let $\mathcal{C}(G)$ be the set of closed subsets of $G$ endowed with the Fell topology.
Define the following subsets of $\mathcal{C}(G)$. 
\begin{align*}
\overline{Y}_u:&=\{A \in \mathcal{C}(G): -P+A \subset A, A \neq \emptyset\},\\
\overline{X}_u:&=\{A \in \overline{Y}_u: -P \subset A\}=\{A \in \overline{Y}_u: 0 \in A\},\\
Y_u:&=\{A \in \overline{Y}_u: A \neq G\},\\
X_u:&=Y_u \cap \overline{X}_u,\\
\overline{X}_{u}^{(0)}:&=\{A \in \overline{Y}_u: A \cap Int(P) \neq \emptyset\}=\{A \in \overline{Y}_u: 0 \in Int(A)\},\\
X_{u}^{(0)}:&=\{A \in Y_u: A \cap Int(P) \neq \emptyset\}.
\end{align*}
Then, $\overline{Y}_u$ is a locally compact, second countable, Hausdorff space which is also a $G$-space, where the action of $G$ on $\overline{Y}_u$ is given by
\[
\overline{Y}_u \times G \ni (A,s) \to A+s \in \overline{Y}_u.\]
Moreover, $\overline{X}_u$ is a compact subset of $\overline{Y}_u$, and $\overline{X}_u+P \subset \overline{X}_u$. Observe that $\overline{X}_{u}^{(0)}$ is an open subset of $\overline{Y}_u$ contained in $\overline{X}_u$, and $\overline{X}_u^{(0)}+P \subset \overline{X}_u^{(0)}$. 
Also, $Y_u$ is an open subset of $\overline{Y}_u$, and $X_u^{(0)}$ is an open subset of $Y_u$ that is contained in $X_u$.

\begin{rmrk}
\label{useful remark}
Let us a record a few facts that we keep appealing to. 
\begin{enumerate}
\item[(1)] The sequence $(\overline{X}_u-na_0)_n$ increases to $\overline{Y}_u$. For a proof of this, we refer the reader to Eq. 2.1 of \cite{Sundar_NYJM} (Pages 1515-1516). Also, $(X_u-na_0)_n \nearrow Y_u$. 
\item[(2)] The sequence $(\overline{X}_{u}+na_0)_{n \geq 1}$ decreases to $\{G\}$. It is clear that $(\overline{X}_u+na_0)_n$ is decreasing, and $\{G\} \subset \bigcap_{n \geq 1}(\overline{X}_u+na_0)$. Suppose $A \in \bigcap_{n \geq 1}(\overline{X}_u+na_0)$. Then, $A-na_0 \in \overline{X}_u$ for every $n \geq 1$, i.e. $na_0 \in A$ for every $n \geq 1$. Since $-P+A \subset A$, $-P+na_0 \subset A$. Since $a_0$ is an order unit, $G=\bigcup_{ n \geq 1}(-P+na_0) \subset A$. Consequently, $A=G$. This proves that 
$\bigcap_{ n \geq 1}(\overline{X}_u+na_0)=\{G\}$. Also, 
 the sequence $(X_u+na_0)_{n \geq 1}$ decreases to $\emptyset$.
\item[(3)] Notice that for $a \in Int(P)$, $X_u+a \subset X_{u}^{(0)}$. Consequently, for every $n \geq 1$, $X_u+(n+1)a_0 \subset X_{u}^{(0)}+na_0$, and $X_u-na_0 \subset X_{u}^{(0)}-(n+1)a_0$ for every $n \geq 1$. Therefore, by $(1)$,  and $(2)$,  $(X_u^{(0)}+na_0)_n \searrow \emptyset$, and $(X_u^{(0)}-na_0)_n \nearrow Y_u$. 

\item[(4)] The collection $\{\overline{X}_u+s: s \in G \}$ generates the Borel $\sigma$-algebra of $\overline{Y}_u$, and the collection $\{(\overline{X}_u+s)\cap  \overline{X}_u: s \in G\}$ generates the Borel $\sigma$-algebra of $\overline{X}_u$. Also, the collection $\{X_u+s: s \in G\}$ generates the Borel $\sigma$-algebra of $Y_u$, and the collection $\{(X_u+s)\cap X_u:s \in G\}$ generates the Borel $\sigma$-algebra of $X_u$. This follows from Lemma 2.1 of \cite{Sundar_NYJM}.

\end{enumerate}
\end{rmrk}

Let $V:=\{V_a\}_{a \in P}$ be an isometric representation of $P$ on a Hilbert space $H$ with commuting range projections. For $s \in G$, write $s=a-b$ with $a,b \in P$, and set $W_s:=V_b^*V_a$. Then, by Prop. 3.4 of \cite{Sundar_Ore}, $W_s$ is a well-defined partial isometry. For $s \in G$, set $E_s:=W_sW_{s}^{*}$. Thanks to Prop. 3.4 of \cite{Sundar_Ore}, $\{E_s: s \in G\}$ is a commuting family of projections.

\begin{ppsn}
\label{spectral measure}
Keep the foregoing notation. There exists a unique projection valued measure $R$ on $\overline{X}_u$ that takes values in $B(H)$ such that for a Borel set $E \subset \overline{X}_u$ and $s \in G$,
\[
W_sR(E)W_s^{*}=R((E+s)\cap \overline{X}_u).\]
\end{ppsn}
\textit{Proof.} For the existence of such a projection valued measure, we refer the reader to  Lemma 7.1 of \cite{Sundar_Ore}. Suppose $R_1$ and $R_2$ are two projection valued measures on $\overline{X}_u$ such that 
for $i=1,2$, 
\[
W_sR_i(E)W_s^{*}=R_i((E+s)\cap \overline{X}_u)\]
for every $s \in G$ and for every Borel set $E \subset \overline{X}_u$. Then, for $s \in G$,
\[
R_1((\overline{X}_u+s)\cap \overline{X}_u)=W_sW_{s}^{*}=R_2((\overline{X}_u+s)\cap \overline{X}_u).\]
As observed in Remark \ref{useful remark}, $\{(\overline{X}_u+s)\cap \overline{X}_u: s \in G\}$ generates the Borel $\sigma$-algebra of $\overline{X}_u$. Hence, $R_1=R_2$. This completes the proof. \hfill $\Box$.

For an isometric representation $V:=\{V_a\}_{a \in P}$ of $P$ with commuting range projections, we call the projection valued measure $R$, given by Prop. \ref{spectral measure}, the projection valued measure associated to $V$. 

\begin{ppsn}
Let $V:=\{V_a\}_{a \in P}$ be an isometric representation of $P$ on a Hilbert space $H$ with commuting range projections, and let $R$ be the projection  valued measure associated to $V$. 
Then, $V$ is pure if and only if $R$ is supported on $X_u$, i.e. $R(\{G\})=0$. 
\end{ppsn}
\textit{Proof.}  Since $\{na_0: n \geq 1\}$ is a cofinal sequence,  $\bigcap_{a \in P}Ran(V_a)=\bigcap_{a \in P}Ran(V_{na_0})$. The covariance relation 
\[V_{na_0}R(E)V_{na_0}^{*}=R(E+na_0)\] for every Borel subset $E \subset \overline{X}_u$ implies that the orthogonal projection onto $Ran(V_{na_0})$ is $R(\overline{X}_u+na_0)$. 
By Remark \ref{useful remark}, $(\overline{X}_u+na_0)_{n \geq 1} \searrow \{G\}$. Thus, $R(\overline{X}_u+na_0) \searrow R(\{G\})$. Hence, $\bigcap_{n \geq 1}Ran(V_{na_0})=\{0\}$ if and only if $R(\{G\})=0$. This completes the proof. \hfill $\Box$

Next, we introduce notation to denote various subsets of the set of non-zero, Radon measures on $Y_u$. Denote by $\mathcal{M}(Y_u)$ the set of non-zero, Radon measures on $Y_u$. Let $\beta \in \bbr$.  Set
\begin{align*}
\mathcal{M}_\beta(Y_u):&=\{ \mu \in \mathcal{M}(Y_u): \textrm{$\mu$ is $e^{-\beta c}$-conformal}\},\\
\mathcal{M}_{a_0}(Y_u):&=\{\mu \in \mathcal{M}(Y_u): \mu(X_u\backslash(X_u+a_0))=1\},\\
\mathcal{M}_{\beta,a_0}(Y_u):&=\mathcal{M}_\beta(Y_u) \cap \mathcal{M}_{a_0}(Y_u),\\
\mathcal{M}_{e}(Y_u):&=\{\mu \in \mathcal{M}(Y_u): \textrm{$\mu$ is ergodic for the $G$-action on $Y_u$}\},\\
\mathcal{M}_{e,\beta}(Y_u):&=\mathcal{M}_\beta(Y_u) \cap \mathcal{M}_{e}(Y_u),\\
\mathcal{M}_{e,\beta,a_0}(Y_u):&=\mathcal{M}_{e,\beta}(Y_u) \cap \mathcal{M}_{a_0}(Y_u).
\end{align*}

\begin{lmma}
\label{compact closure of lshape}
For every $a \in P$, $X_u\backslash(X_u+a)$ has compact closure. Suppose $K \subset X_u$ is a compact subset. Then, there exists $a \in P$ such that $K \subset X_u\backslash(X_u+a)$. 
\end{lmma}
\textit{Proof.}  Let $a \in P$ be given. Observe that $X_u\backslash(X_u+a) \subset X_u\backslash(X_u^{(0)}+a)=\overline{X}_u\backslash(\overline{X}_{u}^{(0)}+a)$, and the latter set is compact as $\overline{X}_u$ is compact, and $\overline{X}_u^{(0)}$ is open. Hence, $X_u\backslash(X_u+a)$ has compact closure. 

Let $K \subset X_u$ be a compact subset.  Note that $X_u\backslash(X_u+na_0)$ is open in $X_u$ for every $n \geq 1$, and by Remark \ref{useful remark}, $(X_u\backslash(X_u+na_0))_n \nearrow X_u$. Thus, there exists $n \geq 1$ such that $K \subset X_u\backslash(X_u+na_0)$. This completes the proof. \hfill $\Box$

\begin{lmma}
\label{lshape has positive measure}
Let $\mu \in \mathcal{M}_{\beta}(Y_u)$ be given. Then, $\mu(X_u\backslash(X_u+a_0)) \in (0,\infty)$. Consequently, there exists a unique $\nu \in \mathcal{M}_{\beta,a_0}(Y_u)$ and a unique positive number $r$ such that $\mu=r\nu$. 
\end{lmma}
\textit{Proof.} it follows from  Lemma \ref{compact closure of lshape} that $\mu(X_u\backslash(X_u+a_0))$ is finite. Suppose that $\mu(X_u\backslash(X_u+a_0))=0$. Then, by the conformality condition, $X_u+ka_0\backslash(X_u+(k+1)a_0)$ is a null set for every $k \geq 0$. By Remark \ref{useful remark}, $(X_u+ka_0)_k \searrow \emptyset$. Hence, 
\[
X_u:=\coprod_{k \geq 0}(X_u+ka_0\backslash(X_u+(k+1)a_0)).\]
Thus, $X_u$ is a null set. The conformality condition forces that  $X_u-na_0$ is a null set for every $n \geq 1$. 
Again by Remark \ref{useful remark}, $Y_u=\bigcup_{ n \geq 1}(X_u-na_0)$. Consequently, $\mu(Y_u)=0$ which is a contradiction. Hence the proof. \hfill $\Box$

\begin{rmrk}
\label{extending functions}
Let $\beta_1,\beta_2$ be real numbers, and let $\mu \in \mathcal{M}_{\beta_1}(Y_u)$. Suppose $\phi:X_u \to \bbc$ is a measurable function such that for every $a \in P$,
\[
\phi(A+a)-\phi(A)=\beta_2 c(a)\]
for $\mu$-almost all $A \in X_u$. 

Using the equality $Y_u:=\bigcup_{a \in P}(X_u-a)=\bigcup_{n \geq 1}(X_u-na_0)$, it is not difficult to prove that there exists a measurable function $\widetilde{\phi}:Y_u \to \bbc$ such that $\widetilde{\phi}|_{X_u}=\phi$, and
for every $s \in G$, 
\[
\widetilde{\phi}(A+s)-\widetilde{\phi}(A)=\beta_2c(s)\]
for $\mu$-almost all $A \in Y_u$. We omit the proof as this is elementary. 

Let  $\beta \in \bbr$ be a real number, and let $\mu$ be a non-zero, Radon  measure on $X_u$ such that $\mu(E+a)=e^{-\beta c(a)} \mu(E)$ for every Borel subset $E\subset X_u$ and for every $a \in P$. Using the fact that $Y_u=\bigcup_{a \in P}(X_u-a)=\bigcup_{a \in P}(X_u^{(0)}-a)$,  it is not difficult to prove that there exists a unique non-zero, Radon measure $\widetilde{\mu}$ on $Y_u$  such that $\widetilde{\mu}|_{X_u}=\mu$, and $\widetilde{\mu}$ is $e^{-\beta c}$-conformal. We leave the details to the reader. 
\end{rmrk}

Let $\mu \in \mathcal{M}_{\beta}(Y_u)$ be given. Consider the Koopman representation $U:=\{U_s\}_{s \in G}$ of $G$ on $L^2(Y_u,\mu)$, i.e.
\[
U_sf(A):=e^{\frac{\beta c(s)}{2}}f(A-s).\]
We view $L^2(X_u,\mu)$ as a closed subspace of $L^2(Y_u,\mu)$. The fact that $X_u+P \subset X_u$ implies that  $\{U_a\}_{a \in P}$ leaves $L^2(X_u,\mu)$ invariant. For $a \in P$, set 
$$V_a:=U_{a}|_{L^2(X_u,\mu)}.$$ Then, $V:=\{V_a\}_{a \in P}$ is a strongly, continuous semigroup of isometries on $L^2(X_u,\mu)$ with commuting range projections. To denote the dependence 
of $V$ on $\mu$, we denote $V$ by $V^{\mu}$.

\begin{ppsn}
\label{building block}
Let $\beta \in \bbr$, and let $\mu \in \mathcal{M}_\beta (Y_u)$. Then, $V^{\mu}$ is $e^{-\frac{\beta c}{2}}$-conformal. Moreovoer, $\dim \mathcal{A}_\beta(V)=1$ if and only if 
$\mu$ is ergodic. 
\end{ppsn}
\textit{Proof.} Let $V:=V^{\mu}$. Note that, in this case, the space $\mathcal{L}_V$ can be written as 
\[
\mathcal{L}_V:=\{ \phi: X_u \to \bbc: \textrm{$\phi$ is measurable, and for every $a \in P$, $\int_{X_u\backslash(X_u+a)}|\phi(A)|^{2}d\mu(A)<\infty$}\}.\] 
As usual, we identify two functions if they agree almost everywhere. 
Thanks to Lemma \ref{compact closure of lshape}, $\mathcal{L}_V=L^{2}_{loc}(X_u,\mu)$. The action of $P$, $T:=\{T_a\}_{a \in P}$, on $\mathcal{L}_V$ is then given by 
\[
T_a\phi(A)=e^{\frac{-\beta c(a)}{2}}\phi(A+a).\]
Take $\phi=1_{X_u}$. Then, $T_a\phi=e^{\frac{-\beta c(a)}{2}}\phi$ for every $a \in P$. Thus, the isometric representation $V$ is $e^{-\frac{\beta c}{2}}$-conformal. 

Suppose that $\mu$ is ergodic. Let $\phi \in \mathcal{L}_V$ be such that $T_a\phi=e^{-\frac{\beta c(a)}{2}}\phi$. Then, for every $a \in P$,
\[
\phi(A+a)=\phi(A)\]
for almost all $A \in X_u$. Thanks to Remark \ref{extending functions}, there exists a Borel function $\widetilde{\phi}:Y_u \to \bbc$ that extends $\phi$ and 
for every $s \in G$, 
$
\widetilde{\phi}(A+s)=\widetilde{\phi}(A)$ 
for almost all $A \in Y_u$. The ergodicity of $\mu$ implies that there exists $c \in \bbc$ such that $\widetilde{\phi}=c$ a.e. Thus, $\phi=c1_{X_u}$. Hence, $\mathcal{A}_\beta(V^{\mu})$ is $1$-dimensional
if $\mu$ is ergodic. 

Conversely, suppose that $\mathcal{A}_{\beta}(V)$ is $1$-dimensional. Let $\phi:Y_u \to \bbc$ be a bounded Borel function that is $G$-invariant. Then, $\phi|_{X_u} \in \mathcal{A}_\beta(V)$. 
Since $\dim \mathcal{A}_\beta(V)$ is $1$-dimensional, there exists $c \in \bbc$ such that $\phi|_{X_u}=c1_{X_u}$. Thus, $\phi(A)=c$ for almost all $A \in X_u$. 
Since $Y_u=\bigcup_{n \geq 1}(X_u-na_0)$ and $\phi$ is $G$-invariant, it follows that $\phi(A)=c$ for almost all $A \in Y_u$. Hence, $\mu$ is ergodic. \hfill $\Box$

\begin{rmrk}
\label{spectral measure of V}
Let $\mu \in \mathcal{M}_{\beta}(Y_u)$, and let $V:=V^{\mu}$. Let $M:L^{\infty}(X_u,\mu) \to B(L^2(X_u,\mu))$ be the multiplication representation. Suppose $R$ is the projection valued measure 
associated to $V$, then $R(E)=M(1_E)$ for every Borel set $E \subset X_u$.  

To see this, for $s=a-b \in G$, with $a,b \in P$, set $W_s:=V_b^{*}V_a$. Then, for $s \in G$, 
\begin{equation*}
W_sf(A):=\begin{cases}
 e^{\frac{\beta c(s)}{2}}f(A-s) & \mbox{ if
} A-s \in X_u,\cr
   &\cr
   0 &  \mbox{ if } A-s \notin X_u
         \end{cases}
\end{equation*}
for $f \in L^2(X_u,\mu)$. Routine computations show that for a Borel set $E \subset X_u$, \[W_sM(1_E)W_{s}^{*}=M(1_{(E+s)\cap X_u}).\] It follows from Prop. \ref{spectral measure} that $R(E)=M(1_E)$ for every Borel subset $E \subset X_u$. 
\end{rmrk}

\begin{ppsn}
\label{absolute continuity}
Let $\beta_1,\beta_2 \in \bbr$, and suppose $\mu_1 \in \mathcal{M}_{\beta_1}(Y_u)$ and $\mu_2 \in \mathcal{M}_{\beta_2}(Y_u)$. Then, the following are equivalent.
\begin{enumerate}
\item[(1)] The isometric representations $V^{\mu_1}$ and $V^{\mu_2}$ are unitarily equivalent.
\item[(2)] The measures $\mu_1$ and $\mu_2$ are absolutely continuous w.r.t. each other.
\end{enumerate}
\end{ppsn}
\textit{Proof.} Clearly, $(2) \implies (1)$. Suppose that $(1)$ holds. Then, the projection valued measures associated to $V^{\mu_1}$ and $V^{\mu_2}$ are unitarily equivalent. It follows from Remark \ref{spectral measure of V} that the multiplication representations
of $C_0(X_u)$ on $L^2(X_u,\mu_1)$ and on $L^2(X_u,\mu_2)$ are unitarily equivalent. Hence, $\mu_1|_{X_u}$ and $\mu_2|_{X_u}$ are absolutely continuous. Thanks to the conformality condition of $\mu_1$ and $\mu_2$, it follows that, for every $n \geq 1$, $\mu_1|_{X_u-na_0}$ and $\mu_2|_{X_u-na_0}$ are absolutely continuous w.r.t. each other. As observed in Remark \ref{useful remark}, $(X_u-na_0)_n \nearrow Y_u$. Hence, $\mu_1$ and $\mu_2$ are absolutely continuous. This completes the proof. \hfill $\Box$

\begin{ppsn}
\label{ergodicity and irreducibility}
Let $\beta \in \bbr$, and let $\mu \in \mathcal{M}_{\beta}(Y_u)$. Then, the following are equivalent.
\begin{enumerate}
\item[(1)] The isometric representation $V^{\mu}$ is irreducible.
\item[(2)] The measure $\mu$ is ergodic.
\end{enumerate}
\end{ppsn}
\textit{Proof.} Let $V:=V^{\mu}$, and let $M:L^{\infty}(X_u,\mu) \to B(L^2(X_u,\mu))$ be the multiplication representation. Assume that $V$ is irreducible. Suppose $\mu$ is not ergodic. Let $E$ be a $G$-invariant subset of $Y_u$ such that $E$ is neither null nor co-null. We claim that $E \cap X_u$ and $E^c\cap X_u$ have positive measure. 
Suppose $E \cap X_u$ is a null set. Then, for every $n \geq 1$, \[
\mu(E \cap (X_u-na_0))=\mu((E-na_0)\cap (X_u-na_0))=\mu((E \cap X_u)-na_0)=e^{-\beta c(na_0)}\mu(E \cap X_u)=0.\]
Thus, $E \cap (X_u-na_0)$ has measure zero for every $n$. Thus, $E=\bigcup_{n \geq 1}(E \cap (X_u-na_0))$ has measure zero, which is a contradiction. Therefore, $E \cap X_u$ is of positive measure. Similarly, $E^c \cap X_u$ is of positive measure. 

Therefore, $Q:=M(1_{E \cap X_u})$ is a non-trivial projection. Let $a \in P$ be given.  Calculate as follows to observe that for $\xi \in L^2(X_u,\mu)$ and $A \in X_u$,
\begin{align*}
QV_{a}^{*}\xi(A)&=e^{-\frac{\beta c(a)}{2}}1_{E \cap X_u}(A)\xi(A+a)\\
&=e^{-\frac{\beta c(a)}{2}}1_{E}(A) \xi(A+a)\\
&=e^{-\frac{\beta c(a)}{2}}1_{E-a}(A)1_{X_u-a}(A)\xi(A+a) ~~~(\textrm{since $E$ is $G$-invariant, and $X_u-a \subset X_u$})\\
&=e^{-\frac{\beta c(a)}{2}}1_{E\cap X_u}(A+a)\xi(A+a)\\
&=V_{a}^{*}Q\xi(A).
\end{align*}
Thus, $Q$ commutes with $V_a^*$ for every $a \in P$. The fact that $Q$ is self-adjoint implies that $Q$ commutes with $V_a$ for every $a \in P$. Consequently, $Q$ is a non-trivial projection such that $Q \in \{V_a,V_a^*:a \in P\}^{'}$. This contradicts the assumption that $V$ is irreducible. Hence, if $V$ is irreducible, then $\mu$ is ergodic. 

Conversely, suppose $\mu$ is ergodic. Let $U$ be a unitary that lies in the commutant $\{V_a,V_{a}^{*}:a \in P\}^{'}$. Since $U$ intertwines $V$ and $V$, $U$ intertwines the associated projection valued measures. By Remark \ref{spectral measure of V}, we have 
\[
UM(1_E)U^*=M(1_E)\]
for every Borel set $E \subset X_u$. Thus, $U \in M(L^{\infty}(X_u,\mu))'=M(L^{\infty}(X_u,\mu))$. 

Let $\phi:X_u \to \bbt$ be a Borel function such that $U=M(\phi)$.  Thanks to the equality $V_{a}^{*}M(\phi)V_a=M(\phi)$, for every $a \in P$, we have 
\begin{equation}
\label{almost everywhere constant}
\phi(A+a)=\phi(A)
\end{equation}
for almost all $A \in X_u$.  Apply Remark \ref{extending functions} to  extend $\phi$ to a Borel function on $Y_u$, which we again denote by $\phi$, such that
given $s \in G$, 
\[
\phi(A+s)=\phi(A)\]
for almost all $A \in Y_u$.  As $\mu$ is ergodic, it follows that there exists a scalar $c $ such that $\phi=c$ a.e. Therefore, $U$ is a scalar multiple of  the identity operator. 
Hence, $V$ is irreducible. This completes the proof of the implication $(2) \implies (1)$. \hfill $\Box$

Let $\mathcal{I}(P)$ be the collection (up to unitary equivalence) of pure isometric representations of $P$ with commuting range projections. Define
\begin{align*}
\mathcal{I}_{r}(P):&=\{V \in \mathcal{I}(P): \textrm{$V$ is irreducible}\}\\
\mathcal{I}_{\beta}(P):&=\{V \in \mathcal{I}(P): \textrm{$V$ is $e^{-\frac{\beta c}{2}}$-conformal}\}\\
\mathcal{I}_{\beta,r}(P):&=\mathcal{I}_{\beta}(P) \cap \mathcal{I}_{r}(P).
\end{align*}
 We have shown that for $\beta \in \bbr$ and for an ergodic, $e^{-\beta c}$-conformal, non-zero Radon measure $\mu$ on $Y_u$, $V^{\mu} \in \mathcal{I}_{\beta,r}(P)$.

\begin{ppsn}
\label{one-one of the required map}
The map 
\[
\mathcal{M}_{e,\beta,a_0}(Y_u) \ni \mu \to V^{\mu} \in \mathcal{I}_{\beta,r}(P)\]
is $1$-$1$. 
\end{ppsn}
\textit{Proof.} Let $\mu,\nu \in \mathcal{M}_{e,\beta,a_0}(Y_u)$ be such that $V^{\mu}$ is unitarily equivalent to $V^{\nu}$. By Prop. \ref{absolute continuity}, it follows that $\mu$ and $\nu$ are absolutely continuous w.r.t. each other. 
Let $f:Y_u \to (0,\infty)$ be a Borel function such that $d\mu=fd\nu$. The fact that both $\mu$ and $\nu$ are $e^{-\beta c}$-conformal implies that for every $s \in G$, 
\[
f(A+s)=f(A)\]
for $\nu$-almost all $A \in Y_u$. Since $\nu$ is ergodic, there exists $c \in (0,\infty)$ such that $f=c$ a.e. Thus, $\mu=c\nu$. The equality $\mu(X_u\backslash(X_u+a_0))=\nu(X_u\backslash(X_u+a_0))=1$ implies that $c=1$. 
Hence, $\mu=\nu$. This completes the proof. \hfill $\Box$

The main point of this section is to show that the map 
\[
\mathcal{M}_{e,\beta,a_0}(Y_u) \ni \mu \to V^{\mu} \in \mathcal{I}_{\beta,r}(P)\]
is surjective. As mentioned earlier, the proof is inspired by the techniques of \cite{Wu}.

Let $V=\{V_a\}_{a \in P}$ be a pure, isometric representation of $P$ with commuting range projections on a Hilbert space $H$. For $a \in P$, set $E_a:=V_aV_a^{*}$. Let $\beta \in \bbr$, and suppose that $V$ is $e^{-\frac{\beta c}{2}}$-conformal. The isometric representation 
$V$ is fixed until further mention. 

Let $\mathcal{D}_V:=\bigcup_{a \in P}Ker(V_a^*)$, and let
\[
\mathcal{L}_V:=\{\phi: \mathcal{D}_V \to \bbc: \textrm{$\phi$ is linear, and for every $a \in P$, $\phi|_{Ker(V_a^*)}$ is bounded}\}.\] 
Recall that the action of  $P$, $T:=\{T_a\}_{a \in P}$,  on $\mathcal{L}_V$ is given by 
\[
T_a\phi(\xi)=\phi(V_a\xi)\]
for $\phi \in \mathcal{L}_V$ and $\xi \in \mathcal{D}_V$. 
Choose a non-zero $\phi \in \mathcal{L}_V$ such that for every $a \in P$, $T_a\phi=e^{-\frac{\beta c(a)}{2}}\phi$. For $a \in P$, let $\xi_a \in Ker(V_a^*)$ be such that 
\[
\phi(\xi)=\langle \xi|\xi_a \rangle\]
for $\xi \in Ker(V_a^*)$. 
 Thanks to Remark \ref{coherence}, the family $\{\xi_a\}_{a \in P}$ is a coherent family, i.e. for $a,b \in P$, if $a \leq b$, then 
\begin{equation}
\label{coherence1}E_a^{\perp}\xi_b=\xi_a.\end{equation}

 Let $R$ be the projection valued measure associated to $V$.
 Set 
 \[
 H_0:=\overline{span\{R(E)\xi_a: \textrm{$E \subset X_u$ is Borel, and $a \in P$}\}}.\]

The key theorem is the following.
\begin{thm}
\label{direct summand}
With the foregoing notation, the closed subspace $H_0$ is invariant under $\{V_a,V_a^{*}:a \in P\}$, and  there exists $\mu \in \mathcal{M}_{\beta}(Y_u)$ such that $V^\mu$ is unitarily equivalent to $V|_{H_0}$.
\end{thm}
The construction of the required measure $\mu$ for the above theorem is next explained in a series of steps. 
Let $a, b \in P$ be given, and let $\xi \in Ker(V_b^*)$ be given. Then, $V_a\xi \in Ker(V_{a+b}^{*})$. Calculate as follows to observe that 
\begin{align*}
e^{-\frac{\beta c(a)}{2}}\langle \xi|\xi_b\rangle&=e^{-\frac{\beta c(a)}{2}}\phi(\xi)\\
&=T_a\phi(\xi)\\
&=\phi(V_a\xi)\\
&=\langle V_a\xi|\xi_{a+b}\rangle\\
&=\langle \xi|V_a^*\xi_{a+b}\rangle\\
&=\langle \xi|E_b^{\perp}V_{a}^{*}\xi_{a+b}\rangle.\end{align*}
Thus, for $a,b \in P$, 
\begin{equation}
\label{conformality condition}
E_{b}^{\perp}V_{a}^{*}\xi_{a+b}=e^{-\frac{\beta c(a)}{2}}\xi_b.
\end{equation}

Let $R$ be the projection valued measure associated to $V$. Recall that we have the following equations. For every $a \in P$ and for a Borel subset $E \subset X_u$,
\begin{equation}
\label{covariance1}
V_{a}^{*}R(E)V_{a}=R((E-a)\cap X_u)
\end{equation}
and \begin{equation}
\label{covariance2}
V_{a}R(E)V_{a}^{*}=R(E+a).
\end{equation}

 
 \begin{enumerate}
 \item[(1)] Let $a \in P$, and let $\mu_{a}$ be the measure on $X_u$ defined by the equation
 \[
 \mu_{a}(E):=\langle R(E)\xi_a|\xi_a \rangle\]
 for a Borel subset $E$ of $X_u$. Observe that 
 \begin{align*}
 \mu_{a}(X_u+a)&=\langle R(X_u+a)\xi_a|\xi_a \rangle \\
 &=\langle V_{a}V_{a}^{*}\xi_{a}|\xi_{a} \rangle ~~~(\textrm{by Eq. \ref{covariance2}})\\
 &=0~~ (\textrm{since $\xi_a \in Ker(V_a^{*})$}).
  \end{align*}
  In other words, $\mu_{a}$ is supported on $X_u\backslash (X_u+a)$. Moreover, for a Borel subset $E \subset X_u\backslash(X_u+a)$,
  \[\mu_a(E)=\langle R(E)\xi_a|\xi_a\rangle.\]
  
  \item[(2)] The family of measures $\{\mu_{a}: a \in P\}$ is consistent in the following sense. Suppose $E$ is a Borel subset such that $E \subset X_u\backslash(X_u+a)$ and $E \subset X_u\backslash(X_u+b)$ for $a,b \in P$. Then, $\mu_{a}(E)=\mu_{b}(E)$. 
  To see this, first consider the case where $b\geq a$, i.e. $b-a \in P$. Write $b=a+d$ with $d \in P$. Suppose $E \subset X_u\backslash(X_u+a)$. Then, $E \subset X_u\backslash(X_u+a+d)$. Calculate as follows to observe that 
  \begin{align*}
  \mu_{b}(E)&=\langle R(E)\xi_{a+d}|\xi_{a+d} \rangle \\
  &= \langle R(X_u\backslash(X_u+a))R(E)R(X_u\backslash(X_u+a))\xi_{a+d}|\xi_{a+d} \rangle \\
  &=\langle E_a^{\perp}R(E)E_{a}^{\perp}\xi_{a+d}|\xi_{a+d} \rangle ~~~(\textrm{by Eq. \ref{covariance2}})\\
  &=\langle R(E)E_{a}^{\perp}\xi_{a+d}|E_{a}^{\perp}\xi_{a+d}\rangle\\
  &=\langle R(E)\xi_a|\xi_a\rangle~~~~~~(\textrm{by Eq. \ref{coherence1}})\\
  &=\mu_a(E).
  \end{align*}
     
     Let $a,b \in P$ be given. Suppose $E$ is a Borel subset of $X_u$ such that $E$ is contained in $(X_u \backslash(X_u+a))\cap (X_u \backslash(X_u+b)) \subset X_u\backslash(X_u+a+b)$. The reasoning made so far implies that 
     \[\mu_{a}(E)=\mu_{a+b}(E)=\mu_{b}(E).\] Consequently, the family of measures $\{\mu_{a}: a \in P\}$ is consistent. 
     
   \item[(3)] Using the fact that  $(X_u\backslash(X_u+na_0))_n \nearrow X_u$ and the fact that the family of measures $\{\mu_{na_0}:n \geq 1\}$ is consistent, we see that  there exists a unique measure $\mu$ on $X_u$ such that if $E$ is a Borel subset of $X_u\backslash(X_u+na_0)$ for some $n \geq 1$, then \[\mu(E)=\langle R(E)\xi_{na_0}|\xi_{na_0}\rangle.\] Since $\{na_0: n \geq 1\}$ is cofinal and the family of measures $\{\mu_a: a \in P\}$ is consistent, we see that for $a \in P$ and a Borel subset $E \subset X_u\backslash(X_u+a)$, we have
   \begin{equation}
   \label{defining equation1}
   \mu(E)=\langle R(E)\xi_a|\xi_a\rangle.
   \end{equation}
      In particular, for every $a \in P$, $\mu(X_u\backslash(X_u+a))=\langle E_{a}^{\perp}\xi_a|\xi_a\rangle=\langle \xi_a|\xi_a \rangle < \infty$. It follows from Lemma \ref{compact closure of lshape} that $\mu$ is a Radon measure. Since $\xi_a$ is non-zero for some $a \in P$, \[\mu(X_u\backslash(X_u+a))=\langle \xi_a|\xi_a\rangle>0\] for some $a \in P$.  Thus, $\mu$ is a non-zero, Radon measure. 
   
   \item[(4)]    We claim that  for every Borel subset $E \subset X_u$ and $a \in P$, $\mu(E+a)=e^{-\beta c(a)}\mu(E)$. 
      Let $a \in P$ and let $E \subset X_u$ be a Borel subset.   Since $X_u:=\bigcup_{ n \geq 1}(X_u\backslash(X_u+na_0))$,  it suffices to consider the case when $E \subset X_u\backslash(X_u+b)$ for some $b \in P$. 
   Then, \[E+a \subset X_u+a \backslash(X_u+a+b) \subset X_u\backslash(X_u+a+b).\] Calculate as follows to observe that 
   \begin{align*}
   \mu(E+a)&=\langle R(E+a)\xi_{a+b}|\xi_{a+b} \rangle \\
   &=\langle V_{a}R(E)V_{a}^{*}\xi_{a+b}|\xi_{a+b} \rangle\\
   &=\langle R(E)V_{a}^{*}\xi_{a+b}|V_{a}^{*}\xi_{a+b} \rangle\\
   &=\langle E_b^{\perp}R(E)E_{b}^{\perp}V_{a}^{*}\xi_{a+b}|V_{a}^{*}\xi_{a+b}\rangle~~~(\textrm{since $E_b^{\perp}R(E)=R(E)$})\\
   &=\langle R(E)E_{b}^{\perp}V_{a}^{*}\xi_{a+b}|E_{b}^{\perp}V_{a}^{*}\xi_{a+b} \rangle\\
      &= e^{-\beta c(a)}\langle R(E)\xi_b|\xi_b \rangle ~~~(\textrm{by Eq. \ref{conformality condition}})\\
   &=e^{-\beta c(a)}\mu(E).
   \end{align*}
   This proves the claim. 
   
   Thanks to Remark \ref{extending functions}, there exists a unique $e^{-\beta c}$-conformal, non-zero, Radon measure on $Y_u$ that extends $\mu$. We denote the extension again by $\mu$. 
   \end{enumerate}

We keep the foregoing notation in the following proof. 

\textit{Proof of Thm. \ref{direct summand}.}
For $a \in P$, let $\mathcal{B}_0^{a}$ be the set of  Borel subsets $E$ of $X_u$ that are contained in $X_u\backslash(X_u+a)$. 
    Set $\mathcal{B}_0:=\bigcup_{a \in P}\mathcal{B}_0^{a}$. 
    Note that if $E \in \mathcal{B}_0$, then $E+b \in \mathcal{B}_0$ for every $b \in P$. 
    
    Also, if $E \in \mathcal{B}_0$ and $b \in P$, then
    $(E-b)\cap X_u \in \mathcal{B}_0$. To see this, suppose $E \in \mathcal{B}_0$ and $b \in P$. We can choose $n$ large such that $na_0-b \in P$
    and $E \subset X_{u}\backslash X_u+na_0$. Then, $(E-b)\cap X_u \subset ((X_u-b)\backslash(X_u+na_0-b))\cap X_u \subset X_{u}\backslash(X_u+na_0-b)$. 
    
    Let $a \in P$ be given. Since $\xi_a \in Ker(V_a^*)$, and $R(X_u\backslash(X_u+a))$ is the orthogonal projection onto $Ker(V_a^{*})$, we have 
    $R(X_u\backslash(X_u+a))\xi_a=\xi_a$. Hence, for a Borel subset $E \subset X_u$, 
    \begin{equation}
    \label{reduction to lshape}
    R(E)\xi_a=R(E\cap (X_u\backslash(X_u+a)))\xi_a.
    \end{equation}
        Thus,  \[
    H_0:=\overline{span\{R(E)\xi_a :E \in \mathcal{B}_0, a \in P\}}.\]
    We claim that $H_0$ is invariant under $\{V_a,V_{a}^{*}: a \in P\}$. 
    
    Let $a \in P$, $E \in \mathcal{B}_0$ and $b \in P$ be given. Calculate as follows to observe that 
    \begin{align*}
    V_{b}R(E)\xi_a&=V_{b}R(E)V_{b}^{*}V_{b}\xi_a\\
    &=e^{\frac{\beta c(b)}{2}}V_{b}R(E)V_{b}^{*}V_bE_{a}^{\perp}V_{b}^{*}\xi_{a+b} ~~~(\textrm{by Eq. \ref{conformality condition}})\\
    &=e^{\frac{\beta c(b)}{2}}R(E+b)(E_{a+b}^{\perp}\xi_{a+b}-E_{b}^{\perp}\xi_{a+b})(\textrm{by Eq. \ref{covariance2}})\\
    &=e^{\frac{\beta c(b)}{2}}R(E+b)(\xi_{a+b}-\xi_b).\end{align*}
    Hence, $V_bR(E)\xi_a \in H_0$. 
    
       Now, we show $V_b^{*}R(E)\xi_a \in H_0$. By Eq. \ref{reduction to lshape}, we can assume that $E \subset X_u\backslash(X_u+a)$. 
    Calculate as follows to observe that 
    \begin{align*}
    V_{b}^{*}R(E)\xi_a&=V_{b}^{*}V_{b}V_{b}^{*}R(E)E_{a}^{\perp}\xi_{a+b} ~~~(\textrm{by Eq. \ref{coherence1}})\\
    &=V_{b}^{*}V_{b}V_{b}^{*}R(E)R(X_u\backslash(X_u+a))\xi_{a+b}\\
    &=V_{b}^{*}V_{b}V_{b}^{*}R(E)\xi_{a+b} ~~~(\textrm{since $E \subset X_u\backslash(X_u+a)$})\\
    &=V_{b}^{*}R(E)V_bV_{b}^{*}\xi_{a+b}\\
    &=V_{b}^{*}R(E)(V_bE_{a}^{\perp}V_{b}^{*}+E_{a+b})\xi_{a+b}\\
    &=V_{b}^{*}R(E)V_bE_{a}^{\perp}V_{b}^{*}\xi_{a+b} ~~~(\textrm{since $E_{a+b}\xi_{a+b}=0$})\\
    &=e^{-\frac{\beta c(b)}{2}}R((E-b)\cap X_u)\xi_a ~~~(\textrm{by Eq. \ref{conformality condition}}).
     \end{align*}
    Hence, $H_0$ is invariant under $\{V_{a},V_{a}^{*}:a \in P\}$. Let $W$ be the restriction of the isometric representation $V$ to $H_0$. Then, $W$ is a direct summand of $V$.
    We prove that $W$ is unitarily equivalent to $V^{\mu}$. 
 
Let $E \in \mathcal{B}_0$ be given. We claim that $R(E)\xi_a$ is independent of $a$ as long as $E$ is contained in $ X_{u}\backslash (X_u+a)$. Suppose $E \subset X_{u}\backslash(X_u+a)$ and $E \subset X_{u}\backslash(X_u+b)$ for some $a,b \in P$.
Note that $E \subset X_u\backslash(X_u+a+b)$. Now,
\[
R(E)\xi_a=R(E)E_{a}^{\perp}\xi_{a+b}=R(E)R(X_u\backslash(X_u+a))\xi_{a+b}=R(E)\xi_{a+b}.\]
Similarly, $R(E)\xi_b=R(E)\xi_{a+b}$. Hence,
\begin{equation}
\label{independence}
R(E)\xi_a=R(E)\xi_b
\end{equation}

For $n \geq 1$, let $H_{\mu}^{n}:=L^{2}(X_u\backslash(X_u+na_0))$. We view $H_{\mu}^{n}$ as a closed subspace of $L^2(X_u,\mu)$. Note that $H_{\mu}^{n}$ increases to $L^2(X_u,\mu)$. 
Let $n \geq 1$, and let $E,F$ be Borel subsets of $X_u$ contained in $X_{u}\backslash(X_u+na_0)$. By Eq. \ref{defining equation1}, we have
\[
\langle 1_E|1_F\rangle_{H_\mu}=\langle R(E)\xi_{na_0}|R(F)\xi_{na_0}\rangle.\]
Thus, there exists an isometry $U^{n}:H_{\mu}^{n} \to H_0$ such that 
\[
U^{n}(1_E)=R(E)\xi_{na_0}\]
whenever $E$ is a Borel set contained in $X_u\backslash(X_u+na_0)$. 

Thanks to Eq. \ref{independence}, the isometries $U^{n}$ patch together to define a well defined isometry $U:L^2(X_u,\mu) \to H_0$ such that 
\[
U(1_E)=R(E)\xi_{na_0}\]
whenever $E$ is a Borel set contained in $X_u\backslash(X_u+na_0)$. 

Suppose $E$ is a Borel set contained in $X_u\backslash(X_u+a)$ for some $a \in P$. Choose a large $n$ such that $na_0-a \in P$. Then, $E \subset X_{u}\backslash(X_u+na_0)$ and by Eq. \ref{independence}
\begin{equation*}
U(1_E)=R(E)\xi_a.\end{equation*}
It follows from the above equality and  Eq. \ref{reduction to lshape}, that the range of $U$ is $H_0$. 

We claim that $U$ intertwines $V^{\mu}$ and $W$. Let $E$ be a Borel subset of $X_u$ contained in $X_u\backslash(X_u+b)$ for some $b$. Let $a \in P$ be given. 
Calculate as follows to observe that 
\begin{align*}
UV_{a}^{\mu}(1_E)&=e^{\frac{\beta c(a)}{2}}U(1_{E+a})\\
&=e^{\frac{\beta c(a)}{2}}R(E+a)\xi_{a+b} ~~~(\textrm{~~since $E+a \subset X_u\backslash(X_u+a+b)$})\\
&=e^{\frac{\beta c(a)}{2}}V_{a}R(E)V_{a}^{*}\xi_{a+b}\\
&=e^{\frac{\beta c(a)}{2}}V_{a}R(E)R(X_u\backslash(X_u+b))V_{a}^{*}\xi_{a+b}\\
&=e^{\frac{\beta c(a)}{2}}V_{a}R(E)E_{b}^{\perp}V_{a}^{*}\xi_{a+b}\\
&=V_{a}R(E)\xi_b ~~~(\textrm{by Eq. \ref{conformality condition}})\\
&=V_{a}U(1_E).
\end{align*}
This proves the claim and the proof is complete. \hfill $\Box$

\begin{rmrk}
\label{smallest subspace}
Note that $H_0$ is the smallest closed subspace that contains $\{\xi_a: a \in P\}$ and that is invariant under $\{V_a,V_a^{*}:a \in P\}$.

The case $\beta=0$ deserves special emphasis. Let $\xi:P \to H$ be a map. The map $\xi$ is called an additive cocycle if 
\begin{enumerate}
\item[(1)] for $a \in P$, $\xi_a \in Ker(V_a^*)$, and
\item[(2)] $\xi$ satisfies the cocycle equation, i.e for $a,b \in P$, $\xi_{a+b}=\xi_a+V_a\xi_b$.
\end{enumerate} 
Denote the set of additive cocycles by $\mathcal{A}(V)$. It is clear, from the cocycle equation, that  if $\xi \in \mathcal{A}(V)$, then $\{\xi_a\}_{a \in P}$ is coherent. Let $\xi \in \mathcal{A}(V)$. Let $\phi \in \mathcal{D}_V$ be such that 
\[
\phi(\xi)=\langle \xi|\xi_a \rangle\]
for $\xi \in Ker(V_a^*)$.
 Let $a,b \in P$, and let $\xi \in Ker(V_b^*)$. Then, $V_a\xi \in Ker(V_{a+b}^{*})$, and 
\begin{align*}
\phi(V_a\xi)&=\langle V_a\xi|\xi_{a+b} \rangle\\
&=\langle V_a\xi|\xi_a+V_a\xi_b\rangle\\
&=\langle \xi|\xi_b \rangle\\
&=\phi(\xi).\end{align*}
Hence, $T_a\phi=\phi$ for every $a \in P$. It is not difficult to show that every $T=\{T_a\}_{a \in P}$-invariant element arises this way, and we could identify the invariant elements of $\mathcal{D}_V$ with $\mathcal{A}(V)$.

\end{rmrk}

\begin{thm}
\label{translation}
Let $\beta \in \bbr$. The map 
\[
\mathcal{M}_{e,\beta,a_0}(Y_u) \ni \mu \to V^{\mu} \in \mathcal{I}_{\beta,r}(P)\]
is a bijection. 
\end{thm}
\textit{Proof.} It was observed in Prop. \ref{one-one of the required map} that the prescribed map is $1$-$1$. Let $V \in \mathcal{I}_{\beta,r}(P)$ be given. It follows from Thm. \ref{direct summand} that  there exists $\mu \in \mathcal{M}_{\beta}(Y_u)$ such that 
$V^\mu$ is a direct summand of $V$. Thanks to Lemma \ref{lshape has positive measure}, we can assume that $\mu(X_u\backslash(X_u+a_0))=1$. Since $V$ is irreducible, $V$ and $V^{\mu}$ are unitarily equivalent. The ergodicity of $\mu$ follows from Prop. \ref{ergodicity and irreducibility}. The proof is complete. \hfill $\Box$

\begin{ppsn}
\label{finite type}
Let $\beta \in \bbr$, and let $V \in \mathcal{I}_{\beta,r}(P)$. Suppose $\mu \in \mathcal{M}_{e,\beta}(Y_u)$ is such that $V^{\mu}$ is unitarily equivalent to $V$. 
Then, the following are equivalent.
\begin{enumerate}
\item[(1)] The measure class $[\mu]$ contains a $G$-invariant, non-zero, Radon measure on $Y_u$.
\item[(2)] The isometric representation $V$ is $1$-conformal.
\end{enumerate}
\end{ppsn}
\textit{Proof.} The proof follows from Thm. \ref{translation} and Prop. \ref{absolute continuity}. \hfill $\Box$

Let $H \subset G$ be a closed subgroup. Note that $G$ acts on $G/H$ by translations. By a $P$-space in $G/H$, we mean a non-empty, proper, Borel subset $B \subset G/H$ such that $B+P \subset B$.
Let $B$ be a $P$-space in $G/H$. For $a \in P$, let $V_a$ be the isometry on $L^2(B)$ defined by 
\begin{equation*}
V_af(x):=\begin{cases}
 f(x-a) & \mbox{ if
} x-a \in B,\cr
   &\cr
   0 &  \mbox{ if } x-a \notin B
         \end{cases}
\end{equation*}
for $f \in L^2(B)$. The measure that we consider on $B$ is the Haar measure. Then, $V:=\{V_a\}_{a \in P}$ is a strongly continuous semigroup of isometries with commuting range projections. We denote $V$ by $V^B$, and call $V^B$ the isometric representation associated to $B$. 

The proof of the following proposition is clear, and hence omitted. 
\begin{ppsn}
\label{typeonecriteria}
Let $\beta \in \bbr$, and let $V \in \mathcal{I}_{\beta,r}(P)$ be given. Suppose $\mu \in \mathcal{M}_{e,\beta}(Y_u)$ is such that $V^{\mu}$ is unitarily equivalent to $V$. Assume that $\mu$ is concentrated on an orbit of a point $A \in Y_u$.
Denote the stabiliser of $A$ by $G_A$.  Let $\pi:G \to G/G_A$ be the quotient map, and let $B:=-\pi(A)$. Then, $B$ is a $P$-space in $G/G_A$, and  $V$ is unitarily equivalent to $V^B$.
\end{ppsn}

Let $\beta \in \bbr$, and let $V \in \mathcal{I}_{\beta,r}(P)$. For $t \in \bbr$, let $V^t:=\{V^t_a\}_{a \in P}$ be the semigroup of isometries defined by $V^{t}_a:=e^{itc(a)}V_a$. 

\begin{ppsn}
\label{invariant measure}
Keep the foregoing notation. Suppose $\beta \neq 0$.  Let $\mu \in \mathcal{M}_{e,\beta}(Y_u)$ be such that $V$ is unitarily equivalent to $V^{\mu}$. Then, the following are equivalent.
\begin{enumerate}
\item[(1)] The measure class $[\mu]$ contains a $G$-invariant, $\sigma$-finite measure on $Y_u$.
\item[(2)] For every $t \in \bbr$, $V$ and $V^t$ are unitarily equivalent. 
\item[(3)] There exists $\delta>0$ such that $V$ and $V^t$ are unitarily equivalent for every $t \in [0,\delta]$.
\item[(4)] There exists $\delta>0$ such that for every $t \in [0,\delta]$, $V^t$ is $e^{-\frac{\beta c}{2}}$-conformal. 
\end{enumerate}
\end{ppsn}
\textit{Proof.} We can assume that $V=V^{\mu}$. Let $M:L^{\infty}(X_u,\mu) \to B(L^2(X_u,\mu))$ be the multiplication representation.
Suppose that $[\mu]$ contains a $G$-invariant, $\sigma$-finite measure $\nu$ on $Y_u$. Write $d\nu=e^{\beta g}d\mu$ for some Borel function $g:Y_u \to \bbr$. 
The fact that $\nu$ is $G$-invariant and $\mu$ is $e^{-\beta c}$-conformal implies that for every $s \in G$, 
\[
g(A+s)-g(A)=c(s)\]
for almost all $A \in Y_u$. For $t \in \bbr$, set $U_t:=M(e^{itg})$. Then, it is routine to verify that for $a \in P$, 
\[
U_tV_{a}U_{t}^{*}=e^{itc(a)}V_a.\]
Hence, for every $t \in \bbr$, $V$ and $V^t$ are unitarily equivalent. 

Conversely, assume that $(2)$ holds. Choose a dense, countable subsemigroup $\Gamma_{+}$ of $P$, and let $C^{*}(\Gamma_{+})$ be the universal, unital $C^{*}$-algebra generated by
a semigroup of isometries $\{v_\gamma\}_{\gamma \in \Gamma_{+}}$. By the universal property of $C^{*}(\Gamma_{+})$, there exists a representation $\pi$ of $C^{*}(\Gamma)$ on $L^2(X_u,\mu)$ such that $\pi(v_\gamma)=V_{\gamma}$ for every $\gamma \in \Gamma_{+}$. Again by the universal property of $C^{*}(\Gamma_{+})$, for every $t \in \bbr$, there exists a representation $\pi_t$ of $C^{*}(\Gamma_{+})$
such that $\pi_t(v_\gamma)=e^{itc(\gamma)}V_{\gamma}$. It is clear that the family $\{\pi_t\}_{t \in \bbr}$ moves `measurably`, i.e. for $\xi,\eta \in L^2(X_u,\mu)$, and $x \in C^{*}(\Gamma_{+})$,  
the map 
\[
\bbr \ni t\to \langle \pi_t(x)\xi|\eta \rangle \in \bbc\]
is measurable. 

The given hypothesis implies that for every $t \in \bbr$, $\pi$ and $\pi_t$ are unitarily equivalent. 
By the corollary to Lemma 4.1.4 of \cite{Arveson_invitation}, for every $t \in \bbr$, there exists a unitary operator $U_t$ on $L^2(X_u,\mu)$ such that 
\begin{enumerate}
\item[(1)] for $\xi,\eta \in L^2(X_u,\mu)$, the map $\bbr \ni t \to \langle U_t\xi|\eta \rangle \in \bbc$ is measurable, and
\item[(2)] for every $t$, $U_t\pi(\cdot)U_t^{*}=\pi_t(\cdot)$.
\end{enumerate}

In particular, for every $t \in \bbr$, and for every $\gamma \in \Gamma_{+}$, $U_tV_{\gamma}U_{t}^{*}=e^{itc(\gamma)}V_{\gamma}$. Since $\Gamma_{+}$ is dense in $P$, it follows that for every $t \in \bbr$ and for every $a \in P$, $U_tV_{a}U_{t}^{*}=e^{itc(a)}V_{a}$. 

Let $s,t \in \bbr$. For $a \in P$,  calculate as follows to observe that 
\begin{align*}
U_{s}U_tV_{a}U_{t}^{*}U_s^{*}&=e^{itc(a)}U_sV_{a}U_{s}^{*}\\
&=e^{itc(a)}e^{isc(a)}V_{a}\\
&=U_{s+t}V_{a}U_{s+t}^{*}.
\end{align*}
Thus, $U_{s+t}^{*}U_sU_t \in \{V_a, V_{a}: a \in P\}^{'}$. Since $V$ is irreducible, there exists $\omega(s,t) \in \bbt$ such that 
\[U_sU_t=\omega(s,t)U_{s+t}.\] 
Clearly, $\omega:\bbr \times \bbr \to \bbt$ is measurable, and is a multiplier. As $H^2(\bbr,\bbt)=0$, it follows that there exists a Borel function $f:\bbr \to \bbt$ such that for $s,t \in \bbt$, $\omega(s,t)=\frac{f(s+t)}{f(s)f(t)}$. Replacing $U_t$ by $f(t)U_t$, we can assume that $U_sU_t=U_{s+t}$ for every $s,t \in \bbr$.

Thus, there exists a strongly continuous (as every weakly measurable $1$-parameter group of unitaries is strongly continuous) $1$-parameter group of unitaries $U:=\{U_t\}_{t \in \bbr}$ on $L^2(X_u,\mu)$ such that for every $a \in P$ and for $t \in \bbr$, 
\begin{equation}
\label{intertwiner1}
U_tV_aU_t^{*}=e^{itc(a)}V_a.\end{equation}
Write $U_t=e^{itD}$ with $D$ being the Stone generator. 

Let $R$ be the projection valued measure associated to $V$. It is also clear that $R$ is the projection valued measure associated to $V^t$ for every $t \in \bbr$. Moreover, $R(E)=M(1_E)$ for every Borel subset  $E \subset X_u$. 
Eq. \ref{intertwiner1} implies that for every $t \in \bbr$, $U_tR(E)U_t^{*}=R(E)$. Hence, $D$ is affiliated to $L^{\infty}(X_u,\mu)$. Thus, $D$ is given by a multiplication operator determined by a Borel function $g:X_u \to \bbr$. 
Then, for $t \in \bbr$, $U_t=M(e^{itg})$. 

Calculate as follows to observe that for $a \in P$, $t \in \bbr$, and $\xi \in L^{2}(X_u,\mu)$,
\begin{align*}
e^{itc(a)}e^{itg(A)}\xi(A)&=e^{itc(a)}V_{a}^{*}V_{a}U_t\xi(A)\\
&=V_{a}^{*}U_{t}V_{a}\xi(A) ~~~(\textrm{by Eq. \ref{intertwiner1}})\\
&=e^{itg(A+a)}\xi(A).
\end{align*}
Thus, for every $a \in P$ and for every $t \in \bbr$, $e^{it(c(a)+g(A))}=e^{itg(A+a)}$ for  almost all $A \in X_u$. Hence, for every $a \in P$,
\[
g(A+a)-g(A)=c(a)\]
for almost all $A \in X_u$. By Remark \ref{extending functions}, there exists a Borel function, which we again denote by $g$, $g:Y_u \to \bbr$ such that for every $s \in G$, 
\[
g(A+s)-g(A)=c(s)\]
for almost all $A \in Y_u$. 

Define a measure $\nu$ on $Y_u$ by 
\[
\nu(E):=\int e^{\beta g(A)}1_{E}(A)d\mu(A).\] Then, $\nu$ is a $G$-invariant, $\sigma$-finite measure that is absolutely continuous w.r.t. $\mu$. This completes the proof of the equivalence $(1) \iff (2)$. 

Clearly, $(2) \implies (3)$. Let \[
T:=\{t \in \bbr: \textrm{$V$ and $V^t$ are unitarily equivalent}\}.\]
Note that $T$ is a subgroup of $\bbr$. Therefore, $(3) \implies (2)$. This completes the proof of the equivalence $(2) \iff (3)$.

Clearly, $(3) \implies (4)$ as $V$ is $e^{-\frac{\beta c}{2}}$-conformal. Assume that $(4)$ holds. Let $\delta>0$ be such that for every $t \in [0,\delta]$, $V^t$ is $e^{-\frac{\beta c}{2}}$-conformal.  Let $t \in [0,\delta]$ be given. 
As in Prop. \ref{building block}, we can write $\mathcal{L}_V=L^{2}_{loc}(X_u,\mu)$. The action of $P$, $T=\{T_a\}_{a \in P}$, corresponding to $V$, on $\mathcal{L}_V$ is given by
\[
T_af(A)=e^{-\frac{\beta c(a)}{2}}f(A+a).\]
Similarly, we can write $\mathcal{L}_{V^t}=L^2_{loc}(X_u,\mu)$. Then, the action of $P$ corresponding to $V^t$, $T^t=\{T_a^t\}_{a \in P}$, is given by 
\[
T_a^{t}f(A)=e^{-itc(a)}e^{-\frac{\beta c(a)}{2}}f(A+a).\]
The hypothesis implies that there exists a non-zero $g \in L^2_{loc}(X_u,\mu)$ such that for every $a \in P$, 
\[
T_a^{t}g=e^{-\frac{\beta c(a)}{2}}g.\]
In other words, for every $a \in P$, \begin{equation}
\label{123}
g(A+a)=e^{it c(a)}g(A)
\end{equation} for almost all $A \in X_u$. Then, for every $a \in P$, $|g(A+a)|=|g(A)|$ for almost all $A \in X_u$. Appealing to Remark \ref{extending functions} and using the fact that $\mu$ is ergodic, we see that there exists $c>0$ such that $|g|=c$ a.e. Without loss of generality, we can assume $|g|=1$. Let $U:=M(g)$, where $M:L^{\infty}(X_u,\mu) \to B(L^2(X_u,\mu))$ is the multiplication representation. Then, by routine computation, we see that for every $a \in P$, 
\[
UV_aU^{*}=V^t_a.\]
This completes the proof of the implication $(4) \implies (3)$. \hfill $\Box$


\begin{thm}
\label{clinching theorem}
Let $\beta$ be a non-zero real number. Let $V:=\{V_a\}_{a \in P}$ be a pure, isometric representation of $P$  on a Hilbert space $H$ with commuting range projections. Assume that $V$ is $e^{-\frac{\beta c}{2}}$-conformal, and $\dim \mathcal{A}_\beta(V)=1$. Assume that there exists a strongly continuous group of unitaries $U:=\{U_t\}_{t \in \bbr}$ on $H$ such that for every $a \in P$ and $t \in \bbr$, $U_tV_aU_t^{*}=e^{itc(a)}V_a$. 

Then, there exists $\mu \in \mathcal{M}_{e,\beta}(Y_u)$ such that 
\begin{enumerate}
\item[(1)] the isometric representation $V^{\mu}$ is a direct summand of $V$, and
\item[(2)] the measure class $[\mu]$ contains a $G$-invariant, $\sigma$-finite measure.
\end{enumerate}
\end{thm}
\textit{Proof.} Let $\mathcal{D}_V:=\bigcup_{a \in P}Ker(V_a^*)$. Let $\phi \in \mathcal{L}_V$ be such that $\phi$ is non-zero, and \begin{equation}
\label{eigen123}
T_a\phi=e^{-\frac{\beta c(a)}{2}}\phi
\end{equation} for every $a \in P$. For every $a \in P$, let $\xi_a \in Ker(V_a^*)$ be such that 
\[
\phi(\xi)=\langle \xi|\xi_a \rangle\]
for every $\xi \in Ker(V_a^*)$. 

Denote the smallest closed subspace that contains $\{\xi_a:a \in P\}$ and that is invariant under $\{V_a,V_a^*: a \in P\}$ by $H_0$. Define an isometric representation $W:=\{W_a\}_{a \in P}$ on $H_0$ by setting
$W_a=V_a|_{H_0}$ for $a \in P$. 
Then, by Remark \ref{smallest subspace} and by Thm. \ref{direct summand}, there exists $\mu \in \mathcal{M}_{ \beta}(Y_u)$ such that $V^{\mu}$ is unitarily equivalent to $W$. 

Write $V=W \oplus W^{'}$ for some isometric representation $W^{'}$. Then, 
\[
1=\dim \mathcal{A}_\beta(V)=\dim \mathcal{A}_\beta(W)+\dim \mathcal{A}_{\beta}(W^{'}).\]
Since $W$ is $e^{-\frac{\beta c}{2}}$-conformal, $\dim \mathcal{A}_{\beta}(W) \geq 1$, and the above equality forces that $\mathcal{A}_{\beta}(W)$ is of dimension one.  Thanks to Prop. \ref{building block}, $\mu$ is ergodic.

For $t \in \bbr$ and for $a \in P$, let $V^t_a=e^{itc(a)}V_a$. Note that for $t \in \bbr$, $U_t$ commutes with $E_a:=V_{a}V_{a}^{*}=V^t_aV^{t*}_a$ for every $a \in P$. Thus, for $a \in P$ and for $t \in \bbr$, $U_t$ maps $Ker(V_a^*)$ onto $Ker(V_a^*)$ for every $a \in P$.
Let \[\mathcal{D}_W:=\bigcup_{a \in P} Ker(W_a^*).\] 
For $t \in \bbr$, let $W^t$ be the isometric representation defined by $W^t_a=e^{itc(a)}W_a$. Clearly, $Ker(W_a^{t*})=Ker(W_a^*)$ for every $a \in P$ and for every $t \in \bbr$. Thus, $\mathcal{D}_{W^t}=\mathcal{D}_W$ for every $t \in \bbr$.

Moreover, $\mathcal{D}_W \subset \mathcal{D}_V$. Choose $a_1 \in P$ such that $\xi_{a_1} \neq 0$. Since $U_t\xi_{a_1} \to \xi_{a_1}$ as $t \to 0$, there exists $\delta>0$ such that 
$\langle \xi_{a_1}|U_t\xi_{a_1} \rangle \neq 0$ for every $t \in [0,\delta]$. Fix $t \in [0,\delta]$. Let $\psi:\mathcal{D}_{W^t} \to \bbc$ be defined by 
\[
\psi(\xi)=\phi(U_t^*\xi).\]
Let $b \in P$. Note that for $\xi \in Ker(W^{t*}_b) \subset Ker(V_b^{*})$, $U_t^{*}\xi \in Ker(V_b^*)$. Thus, for $\xi \in Ker(W^{t*}_b)$, 
\[
\psi(\xi)=\langle U_t^{*}\xi|\xi_b\rangle=\langle \xi|U_t\xi_b\rangle.\]
Thus, $\psi$ is bounded on $Ker(W^{t*}_b)$ for every $b \in P$. 

Also, note that $\psi(\xi_{a_1})=\phi(U_t^{*}\xi_{a_1})=\langle U_t^{*}\xi_{a_1}|\xi_{a_1}\rangle=\langle \xi_{a_1}|U_t\xi_{a_1}\rangle \neq 0$. Thus, $\psi$ is a non-zero linear functional on $\mathcal{D}_{W^t}$. 
Let $a,b \in P$, and let $\xi \in Ker(W^{t*}_b)$ be given. Calculate as follows to observe that 
\begin{align*}
\psi(W_{a}^{t}\xi)&=e^{itc(a)}\psi(W_a\xi)\\
&=e^{itc(a)}\phi(U_t^{*}W_a\xi)\\
&=e^{itc(a)}\phi(U_t^*V_a\xi)\\
&=e^{itc(a)}e^{-itc(a)}\phi(V_aU_t^{*}\xi)\\
&=e^{-\frac{\beta c(a)}{2}}\phi(U_t^{*}\xi)~~~(\textrm{by Eq. \ref{eigen123}})\\
&=e^{-\frac{\beta c(a)}{2}}\psi(\xi).
\end{align*}
Therefore, $W^t$ is $e^{-\frac{\beta c}{2}}$-conformal. Hence, for every $t \in [0,\delta]$, $W^t$ is $e^{-\frac{\beta c}{2}}$-conformal. By Prop. \ref{invariant measure}, it follows that the measure class $[\mu]$ contains a $G$-invariant, $\sigma$-finite measure. 
The proof is complete. \hfill $\Box$

Let $V$ be a pure, isometric representation of $P$ with commuting range projections on a Hilbert space $H$. Let $U:=\{U_s\}_{s \in G}$ be the minimal unitary dilation of $V$ on a Hilbert space $K$. For $s \in G$, let $E_s$ be the orthogonal projection 
onto $U_sH$. Then, $\{E_s:s \in G\}$ is a commuting family of projections. Also, for $s,t \in G$, 
\begin{equation}
\label{conjugacy}
U_sE_tU_s^{*}=E_{t+s}.\end{equation}
Let $\mathcal{D}$ be the von Neumann algebra generated by $\{E_s: s \in G\}$. It is clear from Eq. \ref{conjugacy} that $U_t\mathcal{D}U_t^{*}=\mathcal{D}$ for every $t \in G$. This way, the unitary group $U:=\{U_t\}_{t \in G}$ implements an action of $G$
on $\mathcal{D}$. 

\begin{ppsn}
\label{type III lambda}
Keep the foregoing notation. Suppose $V \in \mathcal{I}_{\beta,r}$, and let $\mu \in \mathcal{M}_{e,\beta}(Y_u)$ be such that $V^{\mu}$ is unitarily equivalent to $V$. For $\lambda \in [0,1]$, $\mu$ is of type $III_\lambda$ if and only if $\mathcal{D} \rtimes G$ is a factor of type $III_\lambda$. 
\end{ppsn}
\textit{Proof.} We can assume $V=V^{\mu}$. Let $U:=\{U_s\}_{s \in G}$ be the Koopman representation of $G$ on $L^{2}(Y_u,\mu)$, and let $H:=L^{2}(X_u,\mu)$. By definition, for $a \in P$, $U_a|_{H}=V_a$. Note that, for $n \geq 1$, $U_{na_0}^{*}H=L^{2}(X_u-na_0)$. Since $\bigcup_{n \geq 1}(X_u-na_0)=Y_u$, it follows that $\bigcup_{n \geq 1}U_{na_0}^{*}H$ is dense in $L^{2}(Y_u,\mu)$. Thus, $U$ satisfies the required properties for it to be the minimal unitary dilation of $V$. 

Let $M:L^{\infty}(Y_u,\mu) \to B(L^2(Y_u,\mu))$ be the multiplication representation. Note that for $s \in G$, $M(1_{X_u+s})$ is the orthogonal projection onto $U_sH$. As $\{X_u+s: s \in G\}$ generates the Borel $\sigma$-algebra of $Y_u$, it follows that the von Neumann algebra generated by $\{M(1_{X_u+s}):s \in G\}$ is $L^{\infty}(Y_u,\mu)$. Thus, $\mathcal{D}=L^{\infty}(Y_u,\mu)$. Clearly, the action of $G$ on $\mathcal{D}$ implemented by the unitary group $\{U_s\}_{s \in G}$ coincides with the translation action of $G$ on $L^{\infty}(Y_u,\mu)$. Thus, $\mathcal{D} \rtimes G=L^{\infty}(Y_u)\rtimes G$. The proof follows. \hfill $\Box$

We end this section by describing a useful recipe to build isometric representations with commuting range projections. 
Let $(Y,\mathcal{B})$ be a standard Borel space on which $G$ acts measurably. 
Suppose $X$ is a measurable subset of $Y$ that is $P$-invariant, i.e. $X+P \subset X$. We assume that $\displaystyle Y=\bigcup_{a \in P}(X-a)$. We call such a  pair $(Y,X)$ a $(G,P)$-space. The $(G,P)$-space $(Y,X)$ is said to be \emph{pure} if 
$\displaystyle \bigcap_{a \in P}(X+a)=\emptyset$. 

Let $(Y,X)$ be a pure $(G,P)$-space. Suppose  $\mu$ is a $G$-invariant, $\sigma$-finite measure on $Y$. For $a \in P$, let $V_a$ be the isometry on $L^2(X,\mu)$ defined by 
\begin{equation*}
V_af(x):=\begin{cases}
 f(x-a) & \mbox{ if
} x-a \in X,\cr
   &\cr
   0 &  \mbox{ if } x-a \notin X
         \end{cases}
\end{equation*}
for $f \in L^2(X)$. Then, $V:=\{V_a\}_{a \in P}$ is a strongly continuous semigroup of isometries with commuting range projections. We call $V$ the isometric representation associated to the triple $(Y,X,\mu)$. 
 Let  $M:L^{\infty}(X,\mu) \to B(L^2(X,\mu))$ be the multiplication
representation. Clearly, $V_{a}V_{a}^{*}=M(1_{X+a})$ for $a \in P$. Since $\{na_0: n \geq 1\}$ is cofinal, the intersection $\bigcap_{ n \geq 1}(X+na_0)=\bigcap_{a \in P}(X+a)=\emptyset$. Thus, $V_{na_0}V_{na_0}^{*} \searrow 0$ strongly. In other words, the isometric representation $V$ is pure.

We keep the above notation till the end of this section.

\begin{rmrk}
\label{extending functions1}
Suppose $\beta \in \bbr$. Let $f:X \to \bbc$ be such that for every $a \in P$, \[f(x+a)=e^{-\frac{\beta c(a)}{2}}f(x)\] for almost all $x \in X$. Using the fact that $Y=\bigcup_{a \in P}(X-a)=\bigcup_{a \in P}(X-na_0)$, it is not 
difficult to prove that there exists a measurable function $\widetilde{f}:Y \to \bbc$ such that
\begin{enumerate}
\item[(1)] for every $s \in G$, $\widetilde{f}(y+s)=e^{-\frac{\beta c(s)}{2}}\widetilde{f}(y)$, and
\item[(2)] for $x \in X$, $\widetilde{f}(x)=f(x)$. 
\end{enumerate}
\end{rmrk}

The proof of the next lemma is omitted as it is elementary. 
\begin{lmma}
\label{1 dimension}
Suppose $\mu$ is ergodic, and let $\beta \in \bbr$. For $i=1,2$, let $g_i:Y \to \bbc$ be a non-zero measurable function such that for every $s \in G$, 
$g_i(y+s)=e^{-\frac{\beta c(s)}{2}}g_i(y)$ for almost all $y \in Y$. Then, there exists $c \in \bbc$ such that $g_1=cg_2$. 
\end{lmma}

\begin{ppsn}
\label{estimate}
Suppose $\mu$ is ergodic. We have the following.
\begin{enumerate}
\item[(1)] $V$ is $1$-conformal if and only if for every $a \in P$, $\mu(X\backslash(X+a))<\infty$.
\item[(2)] For every $\beta \in \bbr$, $\dim \mathcal{A}_\beta(V) \in \{0,1\}$.
\end{enumerate}
\end{ppsn}
\textit{Proof.} As in Prop. \ref{building block}, we can write $\mathcal{L}_V$ as 
\[
\mathcal{L}_V:=\{f:X \to \bbc: \textrm{$f$ is measurable, and for every $a \in P$, $\int_{X\backslash(X+a)}|f(x)|^{2}dx<\infty$}\}.\]
The action of $P$, $T:=\{T_a\}_{a \in P}$, on $\mathcal{L}_V$ is then given by 
\[
T_af(x)=f(x+a)\]
for $a \in P$ and $f \in \mathcal{L}_V$. 

Suppose that $V$ is $1$-conformal. Then, there exists a non-zero $f \in \mathcal{L}_V$ such that for every $a \in P$, $f(x+a)=f(x)$ for almost all $x \in X$.
Appealing to Remark \ref{extending functions1} and using the fact that $\mu$ is ergodic, we see that there exists a non-zero complex number $c$ such that 
$f(x)=c$ for almost all $x \in X$. The fact that for every $a \in P$, 
\[
\int_{X\backslash(X+a)}|f(x)|^{2}dx<\infty\]
implies that $\mu(X\backslash(X+a))<\infty$ for every $a \in P$. Conversely, suppose $\mu(X\backslash(X+a))<\infty$ for every $a \in P$. Then, $1_{X} \in \mathcal{L}_V$, and
$T_a(1_X)=1_X$. Hence, $V$ is $1$-conformal. This completes the proof of $(1)$. 

Let $\beta \in \bbr$ be given. Suppose that $\mathcal{A}_{\beta}(V)$ is non-zero. Let $g_1,g_2 \in \mathcal{A}_{\beta}(V)$ be two non-zero elements. Then, for $i=1,2$, and for every $a \in P$,
$g_i(x+a)=e^{-\frac{\beta c(a)}{2}}g_i(x)$ for almost all $x \in X$. Appealing to Remark \ref{extending functions1} and to Lemma \ref{1 dimension}, we see that there exists a non-zero complex number 
such that $g_1=cg_2$. Hence, $\dim \mathcal{A}_{\beta}(V)=1$. \hfill $\Box$

\begin{ppsn}
\label{minimal unitary dilation}
The Koopman representation of $G$ on $L^2(Y,\mu)$ is the minimal unitary dilation of $V$. 
\end{ppsn}
\textit{Proof.} Let $U:=\{U_s\}_{s \in G}$ be the Koopman representation of $G$ on $L^2(Y_u,\mu)$.  Recall that for $s \in G$ and for $f \in L^2(Y_u,\mu)$, 
\[
U_sf(y):=f(y-s).\]
View $L^2(X,\mu)$ as a subspace of $L^2(Y,\mu)$ in the obvious way.
Clearly, for $a \in P$ and $f \in L^2(X,\mu)$, $U_a(f)=V_a(f)$. 
Note that, for $n \geq 1$, $U_{na_0}^{*}L^2(X,\mu)=L^{2}(X-na_0,\mu)$. Since $\{na_0:n \geq 1\}$ is cofinal, $(X-na_0) \nearrow Y$. Consequently, we have that the union $$\displaystyle \bigcup_{ n \geq 1}U_{na_0}^{*}L^2(X,\mu)=\bigcup_{ n \geq 1}L^{2}(X-na_0,\mu)$$ is dense in $L^2(Y,\mu)$. 
Thus, $U=\{U_s\}_{s \in G}$ satisfies all the required properties for it to be a minimal unitary dilation of $V$. Since the minimal unitary dilation is unique, up to unitary equivalence, $U:=\{U_s\}_{s \in G}$ is the minimal unitary dilation of $V$. \hfill $\Box$

\section{Type $II_\infty$ examples}

 Thanks to Prop. \ref{symmetry}, it suffices to prove Thm. \ref{main} for $\beta>0$. Let $\beta$ be a positive real number
fixed for the rest of this section.  For the rest of the paper, let $G:=\bbz^2$, and let $P:=\bbn^2$. Fix a positive real number $\theta$, and let $c:\bbz^2 \to \bbr$ be the homomorphism defined by $c(m,n)=m+n\theta$. 
 Recall that the potential $\chi:\Omega \to \bbr$ is defined by 
\begin{equation*}
\chi(x):=\begin{cases}
 1 & \mbox{ if
} x_{-1} =0,\cr
   &\cr
   - \theta &  \mbox{ if } x_{-1}=1.
         \end{cases}
\end{equation*}

Set $e_1:=(1,0)$, $e_2:=(0,1)$, $v_1:=e_1$, and $v_2:=e_1+e_2$. Recall that 
\begin{align*}
Y_u&:=\{A \subset \bbz^2:-\bbn^2+A \subset A, A \neq \emptyset, A \neq \bbz^2\},\\
X_u&:=\{A \in Y_u: 0 \in A\}.
\end{align*}

Suppose $\theta=\frac{p}{q}$ is rational, where $gcd(p,q)=1$. Choose a matrix $\begin{bmatrix}
x & y \\
z & w\\
\end{bmatrix} \in SL_2(\bbz)$ such that $x,y,z,w \geq 0$, and 
\[
x+z=q;~y+w=p.\]
Let $\phi:\bbz^2 \to \bbz^2$ be the isomorphism that corresponds to the matrix $\begin{bmatrix}
x & y \\
z & w
\end{bmatrix}$. Note that $\phi(\bbn^2) \subset \bbn^2$. We denote the homomorphism $\bbz^2 \to \bbr$ that sends $e_1 \to 1$ and $e_2 \to \theta$ by $c_{\theta}$. If $\theta=1$, we denote $c_\theta$ simply by $c$. It is clear that $\frac{1}{q}( c \circ \phi)=c_\theta$.

Let $Y:=Y_u$, and define a $\bbz^2$-action on $Y$ by \[
A\oplus (m,n)=A+\phi(m,n).\]

\begin{ppsn}
\label{reduction to 1}
With the above notation, we have the following. 
\begin{enumerate}
\item[(1)] If $\mu$ is an $e^{-\frac{\beta c}{q}}$-conformal measure on $Y_u$, then $\mu$ is an $e^{-\beta c_\theta}$-conformal measure on the $\bbz^2$-space $Y$. 
\item[(2)] Let $\Phi:Y \to Y_u$ be the map defined by \[\Phi(A)=\phi^{-1}(A).\]
Then, the map $\Phi$ is a $\bbz^2$-equivariant, embedding which is proper, i.e. $\Phi^{-1}(K)$ is compact for every compact subset $K \subset Y_u$. Consequently, if $\mu$ is an ergodic,  $e^{-\frac{\beta c}{q}}$-conformal, non-zero, Radon measure on $Y_u$ , then $\mu \circ \Phi^{-1}$ is an ergodic, $e^{-\beta c_\theta}$-conformal, non-zero, Radon measure on $Y_u$. Moreover, $\mu$ and $\mu \circ \Phi^{-1}$ have the same Krieger type. 
 \end{enumerate}
\end{ppsn}
\textit{Proof.} The only thing that requires verification is the fact that $\Phi$ is proper. Let $K \subset Y_u$ be a compact subset. Since $(X_u-nv_2)_n$ is an increasing sequence of open sets which increases to $Y_u$, it follows that there exists $n \geq 1$ such that $K \subset X_u-nv_2$. Note that 
\[
\Phi^{-1}(K)=\Phi^{-1}(K+nv_2)\ominus nv_2.\]
Thus, it suffices to prove that $\Phi^{-1}(K+nv_2)$ is compact. In other words, we can assume that $K \subset X_u$. 

By Lemma \ref{compact closure of lshape}, there exists $(m,n) \in \bbn^2$ such that $K \subset X_u\backslash(X_u+(m,n))$. Hence, it suffices to show that $\Phi^{-1}(X_u\backslash(X_u+(m,n)))$ has compact closure. Observe that $\Phi^{-1}(X_u)=X_u$. Thus,
$\Phi^{-1}(X_u\backslash(X_u+(m,n))=X_u\backslash(X_u\oplus (m,n))=X_u\backslash(X_u+\phi(m,n))$. Since $\phi(m,n) \in \bbn^2$, it follows from Lemma \ref{compact closure of lshape} that $X_u\backslash(X_u+\phi(m,n))$ has compact closure. 
Hence, $\Phi^{-1}(X_u\backslash(X_u+(m,n)))$ has compact closure.   \hfill $\Box$

\begin{rmrk}
\label{theta=1}
In view of Prop. \ref{reduction to 1} and Prop. \ref{bijection between conformal measures},  Thm. \ref{main} holds for a rational $\theta$ if we establish
its validity for $\theta=1$.

\end{rmrk}

Let $\psi:\bbz^2 \to \bbz^2$ be the isomorphism that corresponds to the matrix $\begin{bmatrix}
0 & 1\\
1 & 0
\end{bmatrix}$. Then, $\psi(\bbn^2)=\bbn^2$. 
The proof of the following proposition is straightforward, and hence omitted. 
\begin{ppsn}
\label{reduction to greater than 1}
With the foregoing notation, the map $\Psi:Y_u \to Y_u$ defined by \[\Psi(A)=\psi^{-1}(A)\] is a homeomorphism. Let $\theta>0$, and let $\delta=\frac{1}{\theta}$. Suppose $\mu$ is an ergodic,  $e^{-\frac{\beta c_\theta}{\theta}}$-conformal, non-zero, Radon measure on $Y_u$. Then, $\mu \circ \Psi $ is an ergodic,  $e^{-\beta c_{\delta}}$-conformal, non-zero, Radon measure on $Y_u$. Moreover, $\mu$ and $\mu \circ \Psi$ have the same Krieger type. 
\end{ppsn}

\begin{rmrk}
\label{theta greater than 1}
In view of Prop. \ref{reduction to greater than 1} and Prop. \ref{bijection between conformal measures}, it suffices to prove  Thm. \ref{main} when $\theta \geq 1$.

\end{rmrk}

\begin{ppsn}
\label{stabiliser is zero}
Let $\mu \in \mathcal{M}_{e,\beta}$. Suppose that $\mu$ is concentrated on an orbit of a point, say $A$. Let $G_A$ be the stabiliser of $A$, i.e.
\[
G_A:=\{(m,n) \in \bbz^2:A+(m,n)=A\}.\]
If $G_A \neq 0$, then $V^{\mu}$ is $1$-conformal.
\end{ppsn}
\textit{Proof.} Lemma \ref{lshape has positive measure} allows us to  assume $\mu(X_u\backslash(X_u+v_2))=1$. As in Prop. \ref{parametrisation of $Y_u$}, we identify $Y_u$ with $\Omega \times \bbz$ via the $\bbz^2$-equivariant homeomorphism 
\[
\Omega \times \bbz \ni (x,t) \to A(x,t) \in Y_u.\]
Recall that the action of $\bbz^2$ on $\Omega\times \bbz$ is given by 
 \begin{align*}
(x,t)+v_1:&=(\tau(x),x_{-1}+t),\\
(x,t)+v_2:&=(x,t+1).
\end{align*}
Thanks to Prop. \ref{bijection between conformal measures}, there exists a unique $e^{-\beta \chi}$-conformal, probability measure $m$ on $\Omega$ such that \[\mu(E \times \{n\})=e^{-\beta n(1+\theta)}m(E)\]
for every Borel set $E \subset  \Omega$ and for every $n \in \bbz$. 

By translating $A$ if necessary, we can assume that  $A=A(x,0)$ for some $x \in \Omega$. The fact that $G_A$ is non-zero implies that $x$ is a periodic point. Let $p$ be the period of $x$, and let $\Omega_m:=\{x,\tau(x),\cdots,\tau^{p-1}(x)\}$. Then, $m$ is concentrated
on $\Omega_m$, and $\mu$ is concentrated on $\Omega_m \times \bbz$. 

Let $dm_0$ be the counting measure
on $\Omega_m$, and let $dn$ be the counting measure on $\bbz$. Then, $dm_0\times dn$ is a $\bbz^2$-invariant, non-zero, Radon measure on $\Omega_m \times \bbz$, and hence on $\Omega \times \bbz$ (as $\Omega_m$ is a compact subset of $\Omega$)
that is equivalent to $\mu$. It follows from Prop. \ref{finite type} that $V^{\mu}$ is $1$-conformal. This completes the proof. \hfill $\Box$

  To proceed further, we need the following two facts from ergodic theory.
\begin{enumerate}
\item[(A1)] There exists an infinite measure preserving, conservative, ergodic transformation $T$ such that the maximal spectral type of the Koopman operator $U_T$ and the Lebesgue measure are mutulally  singular.
\item[(A2)] Given an irrational $\alpha$, there exists an infinite measure preserving, conservative, ergodic transformation $T$ such that $T \times R_\alpha$ is ergodic and the maximal spectral type of $U_{T \times R_\alpha}$
 and the Lebesgue measure are mutually singular. Here, $R_\alpha$ is the irrational rotation on $[0,1)$ defined by $R_\alpha(x)=x+\alpha \mod 1$, and the measure on $[0,1)$ is the Lebesgue measure. 
\end{enumerate}
It follows from the results of \cite{Host} and \cite{Aaronson_Nadkarni} that Kakutani towers over adding machine provide such transformations. We make a small digression into these issues.

Let \[\mathcal{D}:=\{(x_1,x_2,\cdots): x_n \in \{0,1\}\}\]  be the group of dyadic integers, and let $S$ be the transformation on $\mathcal{D}$ that corresponds to addition by $1$. 
The measure that we consider on $\mathcal{D}$ is the Haar measure which we denote by $d\nu$.
Then, $S$ is measure preserving, ergodic and conservative. 
 Remove from $\mathcal{D}$ the sequences that are eventually constant, and denote the resulting set again by $\mathcal{D}$. 

Let $(n_k)_{k=1}^{\infty}$ be a sequence of positive integers such that $n_{k+1}>3n_k$. Define a `height function' $h:\mathcal{D} \to \bbn$ as follows: for $x \in \mathcal{D}$, let $k$ be the least integer such that $x_k=0$. 
Set \[
h(x):=n_k-\sum_{j=1}^{k-1}n_j.\]
Let \[
\mathcal{D}^{h}:=\{(x,n) \in \mathcal{D} \times \bbn: 1 \leq n \leq h(x)\}.\]
Let $T:\mathcal{D}^{h} \to \mathcal{D}^{h}$ be defined by 
\begin{equation*}
T(x,n):=\begin{cases}
 (x,n+1)  & \mbox{ if
} n  \leq h(x)-1,\cr
   &\cr
   (Sx,1) &  \mbox{ if } n=h(x).
         \end{cases}
\end{equation*}
Consider the measure $d\nu dn$, where $dn$ is the counting measure on $\bbn$, on $\mathcal{D}^{h}$. Then, $T$ is a measure preserving, ergodic, conservative transformation on $\mathcal{D}^{h}$. It can be checked that $\mathcal{D}^{h}$ has infinite measure (see  Remark (Pages 67-68) in \cite{Host}). Let $U_T$ be the Koopman operator  associated to $T$. Recall that $U_T$ is the unitary operator on $L^2(\mathcal{D}^h)$ defined by
\[
U_T(f)=f \circ T^{-1}.\]
 In \cite{Host}, it was shown (Thm. 4.4.1, Page 66 of \cite{Host})  that the maximal spectral type of $U_T$ is the Riesz product 
\[
\prod_{k=1}^{\infty}(1+\cos(2\pi n_kt)).\]
It is well known (\cite{Peyriere}) that the Riesz product $\displaystyle \prod_{k=1}^{\infty}(1+\cos(2\pi n_kt))$ and the Lebesgue measure are mutually singular. This ensures (A1). 

To ensure (A2), the sequence $(n_k)$ has to be chosen appropriately. Let $\alpha$ be an irrational number. Choose a sequence $(n_k)$ of positive integers such that 
$n_{k+1}>3n_k$ and $e^{2\pi i n_k \alpha} \to e^{2 \pi i \alpha}$ as $k \to \infty$. Thanks to Thm. 2.2 of \cite{Aaronson_Nadkarni}, the set of $L^{\infty}$-eigen values of $T$, denoted $e(T)$, is given 
by 
\begin{equation}
\label{eigen values}
e(T):=\{e^{2\pi i s} \in \bbt: \sum_{k=1}^{\infty}|1-e^{2\pi i n_k s}|^{2}<\infty\}.\end{equation}
The fact that $e^{2 \pi i n_k \alpha} \to e^{2 \pi i \alpha}$ implies that $e^{2 \pi i n_k \ell \alpha} \to e^{2 \pi i \ell \alpha}$ for every $\ell \in \bbz \backslash \{0\}$. Since $\alpha$ is irrational, it is clear from Eq. \ref{eigen values} that $e^{2 \pi i \ell \alpha} \notin e(T)$ for every $\ell \in \bbz \backslash\{0\}$. By the  ergodic multiplier theorem (Thm. 2.7.1 of \cite{Aaronson}), it follows that $T \times R_\alpha$ is ergodic. 

Note that \[
U_{T \times R_\alpha}=\bigoplus_{\ell \in \bbz} e^{2 \pi i \ell \alpha}U_T.\]
It follows from the above equality and the fact that the maximal spectral of $U_T$ and the Lebesgue measure are singular that the maximal spectral type of $U_{T \times R_\alpha}$ and the Lebesgue measure are singular. This ensures (A2).

Now, we proceed towards proving Thm. \ref{main}. 

\textbf{The case $\theta=1$:} Let $\beta>0$ be a fixed real number. Let $c:\bbz^2 \to \bbr$ be the homomorphism defined by $c(m,n)=m+n$.

Fix a $\sigma$-finite measure space $(X,\mathcal{B},\nu)$ and an invertible, measure preserving transformation $T:X \to X$ such that 
\begin{enumerate}
\item[(1)] $\nu(X)=\infty$, and
\item[(2)] $T$ is ergodic, conservative and the maximal spectral type of the Koopman operator $U_T$ and the Lebesgue measure on $\bbt$ are singular.
\end{enumerate}

Fix a measurable set $A \subset X$ such that $0 < \nu (A) <\infty$. Since $T$ is conservative, $\bigcup_{n=0}^{\infty}T^n(A)$ and $X$ are equal up to a set of measure zero. Set $E_0:=A$, and for $n \geq 1$, set $E_{n}:=T^{n}A\backslash(A \cup TA \cup \cdots \cup T^{n-1}A)$. Then, $\{E_n\}_{n \geq 0}$ is a collection of disjoint measurable subsets. Let $\displaystyle E:=\coprod_{n=0}^{\infty}E_n$, and set $F:=E^{c}$. 

Let $\displaystyle f:=\sum_{n=1}^{\infty}-n1_{E_n}$. Since $T$ is measure preserving, $\nu(E_n) \leq \nu(A)$ for every $n \geq 0$. Consequently, 
\begin{equation}
\label{integrability}
\int e^{\beta_0 f}d\nu < \infty
\end{equation}
for every $\beta_0>0$.

Let $x \in X$. Suppose $x \in E$. Choose $n$ such that $x \in E_n$. Then, $\displaystyle Tx \in \coprod_{m=0}^{n+1}E_m$. Clearly, $f(Tx)-f(x) \geq -1$. 
Suppose $x \in F$. Then, $Tx \in F \cup E_0$. In this case, $f(Tx)-f(x)=0$. Thus, \[f(Tx)-f(x) \geq -1\]for every $x \in X$.

Let $\overline{Y}:=X \times \bbz$. Consider the measure $d\nu dn$ on $\overline{Y}$, where $dn$ is the counting measure on $\bbz$. Let 
\[
\overline{X}:=\{(x,n) \in X \times \bbz: n \geq -f(x)\}.\]
Define a $\bbz^2$-action on $\overline{Y}$ by 
\begin{align*}
(x,n)+e_1:&=(x,n+1),\\
(x,n)+e_2:&=(Tx,n+1).
\end{align*}
Clearly, the $\bbz^2$-action on $\overline{Y}$ is measure preserving. The ergodicity of $T$ implies that the $\bbz^2$-action is ergodic. The fact that $f \circ T - f \geq -1$ implies that $\overline{X}+\bbn^2 \subset \overline{X}$. Note that $(\overline{Y},\overline{X})$ is a pure $(\bbz^2,\bbn^2)$-space. 

Let $V:=\{V_{(m,n)}\}_{(m,n) \in \bbn^2}$ be the isometric representation associated to the triple $(\overline{Y},\overline{X},d\nu dn)$. Recall that $V_{(1,0)}$ and $V_{(0,1)}$ are defined by 
\begin{equation*}
V_{(1,0)}f(x,n):=\begin{cases}
 f(x,n-1)  & \mbox{ if
} n -1 \geq -f(x)\cr
   &\cr
   0 &  \mbox{ if } n-1 <-f(x),
         \end{cases}
\end{equation*}
and
\begin{equation*}
V_{(0,1)}f(x,n):=\begin{cases}
 f(T^{-1}x,n-1)  & \mbox{ if
} n -1 \geq -f(T^{-1}x)\cr
   &\cr
   0 &  \mbox{ if } n-1 <-f(T^{-1}x).
         \end{cases}
\end{equation*}

Let $\xi:\overline{X} \to \bbc$ be defined by $\xi(x,n)=e^{-\frac{\beta n}{2}}$. Note that
\[
\int |\xi(x,n)|^{2}d\nu(x)dn=\int_{X}\Big(\sum_{k \geq -f(x)}e^{-\beta k}\Big)d\nu=\frac{1}{1-e^{-\beta}} \int_{X} e^{ \beta f(x)}d\nu(x)<\infty.\]
Hence, $\xi \in L^2(\overline{X})$. Clearly, \[V_{(1,0)}^{*}\xi=e^{-\frac{\beta}{2}}\xi, ~~
\textrm{and}~~ V_{(0,1)}^{*}\xi=e^{-\frac{\beta}{2}}\xi.\] Hence, for $(m,n) \in \bbn^2$, $V_{(m,n)}^{*}\xi=e^{\frac{-\beta c(m,n)}{2} }\xi$. Thanks to Lemma \ref{eigen vector implies conformality}, $V$ is $e^{-\frac{\beta c}{2}}$-conformal. Also, it follows from Prop. \ref{estimate} that $\dim \mathcal{A}_\beta(V)=1$. 

Note that 
\[
\overline{X}\backslash (\overline{X}+(1,0))=\{(x,n) \in X \times \bbz: -f(x) \leq n < -f(x)+1\}=\{(x,-f(x)):x \in X\}\]
which has infinite measure. By Prop. \ref{estimate}, it follows that $V$ is not $1$-conformal. 

For $t \in \bbr$, let $U_t:L^{2}(\overline{X}) \to L^2(\overline{X})$ be the unitary defined by 
\[
U_tf(x,n)=e^{ it \beta n}f(x,n).\]
Note that $U:=\{U_t\}_{t \in \bbr}$ is a strongly continuous group of unitaries. Also,  for $t \in \bbr$, \[U_tV_{(1,0)}U_{t}^{*}=e^{it\beta}V_{(1,0)}, \textrm{~~and~~}~~ U_tV_{(0,1)}U_t^{*}=e^{it\beta}V_{(0,1)}.\] Hence, for every $(m,n) \in \bbn^2$ and $t \in \bbr$, \[U_tV_{(m,n)}U_t^{*}=e^{itc(m,n)}V_{(m,n)}.\]

Appealing to Thm. \ref{clinching theorem}, we see that there exists an ergodic, $e^{-\beta c}$-conformal, non-zero, Radon measure $\mu$ on $Y_u$ such that 
\begin{enumerate}
\item[(1)] $V^{\mu}$ is a direct summand of $V$, and
\item[(2)] the measure class $[\mu]$ contains a $\sigma$-finite, $\bbz^2$-invariant measure.
\end{enumerate}

\begin{thm}
\label{final1}
The measure $\mu$ is of type $II_\infty$.
\end{thm}
\textit{Proof.} Appealing to Remark \ref{additivity of index} and the fact that $V$ is not $1$-conformal, we conclude that $V^{\mu}$ is not $1$-conformal. By Prop. \ref{finite type}, it follows that the measure class $[\mu]$ does not contain a $\bbz^2$-invariant, Radon measure. 
Thus, $\mu$ is not of type $II_1$. 

Suppose $\mu$ is atomic, and is concentrated on the orbit of a point $A \in Y_u$. Since $V^{\mu}$ is not $1$-conformal, it follows from Prop. \ref{stabiliser is zero} that the stabiliser $G_A$ of $A$ is zero. 
Set $B:=-A$. Then, it follows from Prop. \ref{typeonecriteria} that $V^{\mu}$ is unitarily equivalent to the isometric representation $W$ on $\ell^2(B)$ defined by 
\[
W_{(m,n)}(\delta_{(r,s)})=\delta_{(r+m,s+n)}.\]

Since $W$ is a direct summand of $V$, it follows that the minimal unitary dilation of $W$ is a direct summand of the minimal unitary dilation of $V$. It follows from Prop. \ref{minimal unitary dilation}
that the minimal unitary dilation of $W$ is the regular representation $\lambda:=\{\lambda_{(m,n)}\}_{(m,n) \in \bbz^2}$ of $\bbz^2$ on $\ell^2(\bbz^2)$. Again by Prop. \ref{minimal unitary dilation}, it follows that the minimal unitary dilation of 
$V$ is the Koopman representation $\widetilde{U}:=\{\widetilde{U}_{(m,n)}\}_{(m,n) \in \bbz^2}$ of $\bbz^2$ on $L^2(\overline{Y})$. Hence, $\lambda$ is a subrepresentation of $\widetilde{U}$. Thus, for every $(m,n)$, 
$\widetilde{U}_{(m,n)}$ is of the form 
\[
\widetilde{U}_{(m,n)}=\begin{bmatrix}
                 \lambda_{(m,n)} & 0\\
                 0 & *
                 \end{bmatrix}.\]
Note that $\widetilde{U}_{(-1,1)}=U_T \otimes 1$, where $U_T$ is the Koopman operator associated to the ergodic transformation $T$, and $\lambda_{(-1,1)}$ is unitarily equivalent to $U \otimes 1$, where $U$ is the bilateral shift on $\ell^2(\bbz)$.

Hence, $U_T \otimes 1$ is of the form 
\[
U_{T} \otimes 1=\begin{bmatrix}
         U \otimes 1 & 0\\
         0 & *
         \end{bmatrix}.\]
                  The above equality contradicts the fact the maximal spectral type of $U_T$ and the Lebesgue measure on $\bbt$ are singular. Hence, $\mu$ is not atomic.
                  The proof is complete. \hfill $\Box$

We have now established Thm. \ref{main} when $\theta$ is rational (see Remark \ref{theta=1}). 

Suppose that $\theta$ is irrational. As observed in Prop. \ref{reduction to greater than 1}, we can assume that $\theta>1$. Let $\alpha:=\theta-\lfloor \theta \rfloor$, where $\lfloor \theta \rfloor$ is the integral part of $\theta$.   Let $c:\bbz^2 \to \bbr$ be the homomorphism
defined by $c(m,n)=m+n\theta$.

Let $R_\alpha$ be the map on $[0,1)$ defined by $R_\alpha(x)=x+\alpha \mod 1$. The measure that we consider on $[0,1)$ is the Lebesgue measure.
Fix a $\sigma$-finite measure space $(X,\mathcal{B},\nu)$ and an invertible, measure preserving transformation $T:X \to X$ such that 
\begin{enumerate}
\item[(1)] $\nu(X)=\infty$, and
\item[(2)] $T\times R_\alpha$ is ergodic, conservative and the maximal spectral type of the Koopman operator $U_{T\times R_\alpha}$ and the Lebesgue measure on $\bbt$ are singular.
\end{enumerate}
 
 Fix a measurable set $A \subset X$ such that $0 < \nu (A) <\infty$. Define the sets $(E_n)_{n=0}^{\infty}$ and the function $f:X \to \bbz$ as in the case $\theta=1$. We have 
 $
 f(Tx)-f(x) \geq -1$
 for every $x \in X$. 
 
 Let $\overline{Y}:=X \times \bbr$. Consider the measure $d\nu ds$ on $\overline{Y}$, where $ds$ is the Lebesgue measure on $\bbr$. 
 Set 
 \[
 \overline{X}:=\{(x,s) \in X \times \bbr:s \geq -f(x)\}.\]
 Define a $\bbz^2$-action on $\overline{Y}$ by settting
 \begin{align*}
 (x,s)+ e_1:&=(x,s+1),\\
 (x,s)+e_2:&=(Tx,s+\theta)=(Tx,s+\lfloor \theta \rfloor +\alpha).
 \end{align*}
 The ergodicity of $T \times R_\alpha$ implies that the $\bbz^2$-action on $\overline{Y}$ is ergodic. As in the case $\theta=1$, we have $\overline{X}+\bbn^2 \subset \overline{X}$. Moreover, 
 $(\overline{Y},\overline{X})$ is a pure $(\bbz^2,\bbn^2)$-space. 

Let $V:=\{V_{(m,n)}\}_{(m,n) \in \bbn^2}$ be the isometric representation associated to the triple $(\overline{Y},\overline{X},d\nu ds)$. Recall that $V_{(1,0)}$ and $V_{(0,1)}$ are defined by 
\begin{equation*}
V_{(1,0)}f(x,s):=\begin{cases}
 f(x,s-1)  & \mbox{ if
} s -1 \geq -f(x),\cr
   &\cr
   0 &  \mbox{ if } s-1 <-f(x),
         \end{cases}
\end{equation*}
and
\begin{equation*}
V_{(0,1)}f(x,s):=\begin{cases}
 f(T^{-1}x,s-\theta)  & \mbox{ if
} s -\theta \geq -f(T^{-1}x),\cr
   &\cr
   0 &  \mbox{ if } s-\theta <-f(T^{-1}x).
         \end{cases}
\end{equation*}

Let $\xi:\overline{X} \to \bbc$ be defined by $\xi(x,s)=e^{-\frac{\beta s}{2}}$. It can be verified using Eq. \ref{integrability} that $\xi \in L^2(\overline{X})$. Clearly, \[
V_{(1,0)}^{*}\xi=e^{-\frac{\beta}{2}}\xi,~~\textrm{and~~} V_{(0,1)}^{*}\xi=e^{-\frac{\beta \theta}{2} }\xi.\] It follows from Lemma \ref{eigen vector implies conformality} that $V$ is $e^{-\frac{\beta c}{2}}$-conformal. It follows from Lemma \ref{estimate} that $\mathcal{A}_{\beta}(V)=1$. 

The following facts can be established as in the case $\theta=1$. 
\begin{enumerate}
\item[(a)] The isometric representation $V$ is not $1$-conformal.
\item[(b)] There exists a strongly continuous group of unitaries $U:=\{U_t\}_{t \in \bbr}$ such that for $t \in \bbr$ and $(m,n) \in \bbn^2$, 
\[
U_tV_{(m,n)}U_t^{*}=e^{itc(m,n)}V_{(m,n)}.\]
\end{enumerate}

Appealing to Thm. \ref{clinching theorem}, we see that there exists an ergodic, $e^{-\beta c}$-conformal, non-zero, Radon measure $\mu$ on $Y_u$ such that 
\begin{enumerate}
\item[(1)] $V^{\mu}$ is a direct summand of $V$, and
\item[(2)] the measure class $[\mu]$ contains a $\sigma$-finite, $\bbz^2$-invariant measure.
\end{enumerate}
We will show that $\mu$ is of type $II_\infty$. This requires expressing the isometric representation $V$ in a slightly different form which we do next. 

Let $X_0:=X \times [0,1)$. The measure that we consider on $X_0$ is $d\nu ds$, where $ds$ is the Lebesgue measure on $[0,1)$. Note that the map 
\[
X_0 \times \bbz \ni ((x,s),m) \to (x,s-f(x)+m) \in \overline{Y}\]
is a measure preserving, Borel isomorphism. This way, we identify $\overline{Y}$ with $X_0 \times \bbz$. Then, $\overline{X}$ gets identified with $X_0 \times \bbn$.  Transport the $\bbz^2$-action on $\overline{Y}$  to $X_0 \times \bbz$. Then, the action of $e_1$ is given by
\begin{equation}
\label{e1 action}
((x,s),m)+e_1=((x,s),m+1).\end{equation}
To describe, the action of $e_2$, we introduce some notation. Let $\phi:X \to \bbz$ be defined by \[
\phi(x)=f(Tx)-f(x)+\lfloor \theta \rfloor.\] For $k \geq  0$, let $I_k:=\{x \in X: \phi(x)=k\}$. The action of $e_2$ is given by 
\begin{equation}
\label{e2 action}
((x,s),m)+e_2:=\begin{cases}
 ((Tx,R_\alpha(s)),m+k)  & \mbox{ if
} x \in I_k, ~~\textrm{and~} s \in [0,1-\alpha)\cr
   &\cr
   ((Tx,R_\alpha(s)),m+k+1) &  \mbox{ if } x \in I_k, ~~\textrm{and~} s \in [1-\alpha,1).
         \end{cases}
\end{equation}
Set $I_{-1}=\emptyset$. For $k \geq 0$, let $P_k$ be the multiplication operator on $L^2(X_0 )$ that corresponds to the characteristic function $1_{(I_k \times [0,1-\alpha))\cup (I_{k-1}\times [1-\alpha,1))}$. 
With the Hilbert space $L^2(\overline{X})$  identified with $L^{2}(X_0 \times \bbn) \cong L^{2}(X_0)\otimes \ell^2(\bbn)$, Eq. \ref{e1 action} and Eq. \ref{e2 action} allow us to  write the operators $V_{(1,0)}$ and $V_{(0,1)}$ as 
\begin{equation}
\label{operator1}
V_{(1,0)}=1\otimes S,
\end{equation}
and 
\begin{equation}
\label{operator2}
V_{(0,1)}=\sum_{k=0}^{\infty}U_{T\times R_\alpha}P_k\otimes S^k.
\end{equation}
In the above equations, $S$ stands for the unilateral shift on $\ell^2(\bbn)$.

\begin{thm}
\label{final2}
The measure $\mu$ is of type $II_\infty$. 
\end{thm}
\textit{Proof.} Suppose $\mu$ is not of type $II_\infty$. Then, arguing as in Prop. \ref{final1}, we see that $\mu$ is supported on an orbit of a point, say $A$, whose stabiliser is zero.
Set $B:=-A$. Then, $V^{\mu}$ is unitarily equivalent to the isometric representation  $W:=\{(W_{(m,n)}\}_{(m,n)\in \bbz^2}$  on $\ell^2(B)$ defined by 
\[
W_{(m,n)}(\delta_{(r,s)})=\delta_{(r+m,n+s)}.\]

Since $W$ is a direct summand of $V$ and $V_{(1,0)}$ is a pure isometry, $W_{(1,0)}$ is a pure isometry. This implies that
\begin{equation}
 \label{purity}
\bigcap_{m \geq 1}(B+(m,0))=\emptyset.
\end{equation}

The fact that $B+e_1 \subset B$ implies that given $n \in \bbz$, $\{m \in \bbz: (m,n) \in B\}$ is either empty or is an interval in $\bbz$. 
Eq. \ref{purity} implies that for every $n \in \bbz$, the set $\{m \in \bbz: (m,n) \in B\}$ is bounded below. For $n \in \bbz$, set 
\[
a_n:=\inf \{m \in \bbz: (m,n) \in B\}.\]
Note that $a_n$ could be infinity for some $n$. 

We claim that $B=\{(m,n) \in \bbz^2: m \geq a_n\}$. The inclusion $\subset$ is clear.  Suppose $(m,n) \in \bbz^2$ is such that $m \geq a_n$. Then, $a_n$ is finite. By definition, 
$(a_n,n) \in B$. Since $B+(m-a_n,0) \subset B$, it follows that $(m,n) \in B$. This proves the claim. 

The fact that $B+e_2 \subset B$ implies that the sequence $(a_n)$ is decreasing, i.e. for $n \in \bbz$, $a_{n+1} \leq a_n$. 
Let \[
F:=\{n \in \bbz: a_n <\infty\}.\]
Note that $F$ is non-empty as $B$ is non-empty. Also, since $(a_n)$ is decreasing, it follows that $F+1 \subset F$. Hence, $F$ is an interval. 
For $k \geq 0$, let 
\[
F_k:=\{n \in F: a_n-a_{n+1}=k\}.\]
For $k \geq 0$, let $Q_k:\ell^2(F) \to \ell^2(F)$ be the projection onto $\ell^2(F_k)$. We denote both the right shift on $\ell^2(F)$ and $\ell^2(\bbn)$ by the same letter $S$. 

We identify $B$ with $F \times \bbn$ via the map  \[
F \times \bbn \ni (r,s) \to (s+a_r,r) \in B.\]
Once $\ell^2(B)$ is identified with $\ell^2(F \times \bbn)=\ell^2(F)\otimes \bbn$, the operator $W_{(1,0)}$ on $\ell^2(F \times \bbn)$ is given 
\begin{equation}
\label{operator11}
W_{(1,0)}(\delta_{(r,s)})=\delta_{(r,s+1)}.\end{equation}
In other words, $W_{(1,0)}=1\otimes S$. 
The operator $W_{(0,1)}$ on $\ell^2(F \times \bbn)$ is given by 
\[
W_{(0,1)}(\delta_{(r,s)})=\delta_{(r+1,s+k)}~~~ \textrm{if $r \in F_k$}.\]
In other words, 
\begin{equation}
\label{operator22}
W_{(0,1)}=\sum_{k=0}^{\infty}SQ_k \otimes S^k.\end{equation}

Since $W$ is a direct summand of $V$,  it follows there exists a bounded linear operator $J:\ell^2(F)\otimes \ell^2(\bbn) \to L^2(X_0)\otimes \ell^2(\bbn)$, which is also an isometry, such that 
\[
JW_{(m,n)}=V_{(m,n)}J~~~\textrm{and~~~}JW_{(m,n)}^{*}=V_{(m,n)}^{*}J.\]
Thus, $J(1\otimes S)=(1\otimes S)J$, and $J(1\otimes S^{*})=(1\otimes S^{*})J$. Hence, $J$ is of the form $J=J_0\otimes 1$, where $J_0:\ell^2(F) \to L^{2}(X_0)$ is an isometry. 

Thanks to Eq. \ref{operator2}, Eq. \ref{operator22}, and the equality $(J_0\otimes 1)W_{(0,1)}=V_{(0,1)}(J_0\otimes 1)$, we have 
\[
J_0SQ_k=U_{T \times R_\alpha}P_kJ_0\]
for every $k \geq 0$. Since $\sum_k Q_k=1$ and $\sum_{k=1}P_k=1$, we have $J_0S=U_{T \times R_\alpha}J_0$. 

The equality $(J_0 \otimes 1)W_{(0,1)}^{*}=V_{(0,1)}^{*}(J_0\otimes 1)$ leads to the equation $J_0S^*=U_{T \times R_\alpha}^{*}J_0$. In other words, $U_{T \times R_\alpha}$ is of the form
\[
U_{T \times R_\alpha}=\begin{bmatrix}
S & 0 \\
0 & *
\end{bmatrix}.\]
Hence, $S$ is indeed a unitary operator which forces $F=\bbz$. Let us change notation, and denote the bilateral shift on $\ell^2(\bbz)$ by $U$. Then, $U_{T \times R_\alpha}$ is of the form 
\[
U_{T \times R_\alpha}=\begin{bmatrix}
U & 0 \\
0 & *
\end{bmatrix}.\]
The above equality is a contradiction to our assumption that the maximal spectral type of $U_{T \times R_\alpha}$ and the Lebesgue measure on $\bbt$ are singular. Hence, our initial assumption 
is wrong. Consequently, $\mu$ is of type $II_\infty$. The proof is complete. \hfill $\Box$

Theorem \ref{main} is now a consequence of Thm. \ref{final1}, Thm. \ref{final2}, Prop. \ref{bijection between conformal measures}, Remark \ref{theta=1} and Remark \ref{theta greater than 1}.

\section{Type $III$ examples for a rational $\theta$}
In this section, we prove Thm. \ref{main1}. Thanks to Prop. \ref{symmetry}, Prop. \ref{reduction to 1} and Prop. \ref{bijection between conformal measures}, it suffices to prove Thm. \ref{main1}
under the assumption $\beta>0$ and $\theta=1$. Let $\theta=1$. Let $c:\bbz^2 \to \bbr$ be the homomorphism that sends $e_1 \to 1$ and $e_2 \to 1$. 

Let $(Y,X)$ be a pure $(\bbz^2,\bbn^2)$-space. Recall that this means that $Y$ is a standard Borel space on which $\bbz^2$-acts measurably, $X+\bbn^2 \subset X$, $ \bigcup_{(m,n) \in \bbn^2}(X-(m,n))=Y$, and
$\bigcap_{(m,n) \in \bbn^2}(X+(m,n))=\emptyset$. 
Let $\nu$ be $\sigma$-finite measure on $Y$ which is  $e^{-\beta c}$-conformal. Assume that $\nu(X) \in (0,\infty)$. Let $U:=\{U_{(m,n)}\}_{(m,n) \in \bbz^2}$ be the Koopman representation of $\bbz^2$ on $L^2(Y,\nu)$. 
Note that $X+\bbn^2 \subset X$ implies that $L^{2}(X,\nu)$ is invariant under $\{U_{(m,n)}:(m,n) \in \bbn^2\}$. Set $H:=L^{2}(X,\nu)$, and let $V:=\{V_{(m,n)}\}_{(m,n) \in \bbn^2}$ be the isometric representation on $H$ defined by
\[
V_{(m,n)}:=U_{(m,n)}|_{H}.\]
Notice that for $(m,n) \in \bbn^2$, $V_{(m,n)}^{*}(1_X)=e^{-\frac{\beta c(m,n)}{2}}1_X$. Thus, $V$ is $e^{-\frac{\beta c}{2}}$-conformal. Let $M:L^{\infty}(Y,\nu) \to B(L^2(Y,\nu))$ be the multiplication representation.

\begin{ppsn}
\label{trick to produce type III lambda}
Keep the foregoing notation. Assume that $\nu$ is ergodic. The following are equivalent.
\begin{enumerate}
\item[(1)] The isometric representation $V$ is irreducible.
\item[(2)] The set $\{M(1_{X+(m,n)}):(m,n) \in \bbz^2\}$ generates the von Neumann algebra $L^{\infty}(Y,\nu)$.
\end{enumerate}
Suppose that $V$ is irreducible. Let $\mu \in \mathcal{M}_{e,\beta}(Y_u)$ be such that $V^{\mu}$ is unitarily equivalent to $V$. Then, for $\lambda \in [0,1]$, $\mu$ is of type $III_\lambda$ if and only if $\nu$ is of type $III_\lambda$. 
\end{ppsn}
\textit{Proof.} It can be proved as in Prop. \ref{minimal unitary dilation} that $U:=\{U_{(m,n)}\}_{(m,n) \in \bbz^2}$ is the minimal unitary dilation of $V$. Note that for $(m,n) \in \bbz^2$, the orthogonal projection onto $U_{(m,n)}H$ is $M(1_{X+(m,n)})$. 
Let $\mathcal{D}:=\{M(1_{X+(m,n)}):(m,n) \in \bbz^2\}^{''}$. 
Let 
\begin{align*}
N:&=\{V_{(m,n)},V_{(m,n)}^{*}:(m,n) \in \bbn^2\}^{'}\\
M:&=\mathcal{D}^{'} \cap \{U_{(m,n)}:(m,n) \in \bbz^2\}^{'}.
\end{align*}
We leave it to the reader to check that the map $M \ni T \to T|_{H} \in N$ is a $^*$-algebra isomorphism. Thus, $V$ is irreducible if and only if $M=\bbc$. 

Suppose that $(2)$ holds. Then, $\mathcal{D}^{'}=L^{\infty}(Y,\nu)$. The fact that $M=\bbc$ follows from the ergodicity of $\nu$. This completes the proof of the implication $(2) \implies (1)$. 

Conversely, assume that $(1)$ holds.  Then, $M=\bbc$. Denote the $\sigma$-algebra of $Y$ by $\mathcal{B}$, and  let $\mathcal{F}$ be the $\sigma$-algebra generated by $\{X+(m,n):(m,n) \in \bbz^2\}$. Then, $L^{2}(Y_u,\mathcal{F},\nu)$ is a non-zero closed subspace 
of $L^{2}(Y, \mathcal{B}, \nu)$ that is invariant under  $\{U_{(m,n)}:(m,n) \in \bbz^2\}$ and $\mathcal{D}$. The fact that $M=\bbc$ implies that $L^{2}(Y,\mathcal{F},\nu)=L^{2}(Y,\mathcal{B},\nu)$. 

Let $E \in \mathcal{B}$ be given. For $n \geq 1$, let $E_n:=E \cap (X-(n,n))$. Then, $1_{E_n} \nearrow 1_E$. Let $n \geq 1$ be given. Note that $E_n$ has finite measure as $\nu(X)<\infty$ and as $\nu$ is $e^{-\beta c}$-conformal.  Then, $1_{E_n} \in L^{2}(Y,\mathcal{B},\nu)$. Since $L^{2}(Y,\mathcal{B},\nu)=L^{2}(Y,\mathcal{F},\nu)$, there  exists a $\mathcal{F}$-measurable function $f_n$ such that $1_{E_n}=f_n$ a.e. Since $f_n=1_{E_n}$ a.e., $f_n \in L^{\infty}(Y,\mathcal{F},\nu)$. This implies that  $M(1_{E_n})=M(f_n) \in M(L^{\infty}(Y,\mathcal{F},\nu))=\mathcal{D}$. Since $1_{E_n} \nearrow 1_E$, it follows that $M(1_{E}) \in \mathcal{D}$. Consequently, $\mathcal{D}=L^{\infty}(Y,\mathcal{B},\nu)$. This completes the proof of the implication $(1) \implies (2)$. 

Suppose that $V$ is irreducible. Thanks to the equivalence between $(1)$ and $(2)$, we have  $\mathcal{D}=L^{\infty}(Y,\nu)$. The final assertion now follows from Prop. \ref{type III lambda}.

\begin{ppsn}
\label{Base type III}
Let $\beta>0$, and let $\lambda:=e^{-\beta}$. There exists a pure $(\bbz^2,\bbn^2)$-space $(Y,X)$ and an $e^{-\beta c}$-conformal measure $\nu$ on $Y$ such that 
\begin{enumerate}
\item[(1)] $\nu$ is ergodic, 

\item[(2)] the set $\{M(1_{X+(m,n)}):(m,n) \in \bbz^2\}$ generates $L^{\infty}(Y,\nu)$,
\item[(3)] $0<\nu(X)<\infty$, and
\item[(4)] the von Neumann algebra $L^{\infty}(Y)\rtimes \bbz^2$ is a factor of type $III_\lambda$.
\end{enumerate}
Here, $M:L^{\infty}(Y,\nu) \to B(L^{2}(Y,\nu))$ is the multiplication representation. 
\end{ppsn}
\textit{Proof.} Let $\mathcal{D}:=\{(x_1,x_2,\cdots):x_n \in \{0,1\}\}$ be the group of dyadic integers. Remove from $\mathcal{D}$ the eventually constant sequences, and denote
the resulting set again by $\mathcal{D}$. Let $S$ be the transformation on $\mathcal{D}$ that corresponds to addition by $1$.

Let $p \in (0,\frac{1}{2})$ be such that $\frac{1-p}{p}=e^{\beta}$. Let $\displaystyle m:= \otimes_{k=1}^{\infty}p_k$ be the product measure, where $p_k(\{0\})=1-p$, and $p_k(\{1\})
=p$. Then, 
\[
\frac{d(m \circ S)}{dm}=e^{\beta \phi},\]
where $\phi:\mathcal{D} \to \bbz$ is defined by 
\[
\phi(x):=\min \{n \geq 1: x_n=0\}-2.\]

Let $Y:=\mathcal{D} \times \bbz$, and let $X:=\mathcal{D} \times \bbn$. Define a $\bbz^2$-action on $Y$ by setting
\begin{align*}
(x,k)+e_1:&=(x,k+1),\\
(x,k)+e_2:&=(Sx,\phi(x)+k+1).\end{align*}
Since $\phi \geq -1$, it follows that $X+\bbn^2 \subset X$. Moreover, $(Y,X)$ is a pure $(\bbz^2,\mathbb{N}^{2})$-space. 
Let $\nu$ be the measure on $Y$ defined by $d\nu=e^{-\beta k}dm dk$, where $dk$ is the counting measure on $\bbz$. 
Clearly, $\nu$ is $e^{-\beta c}$-conformal, and $\nu(X) \in (0,\infty)$. 

The ergodicity of $\nu$ follows from that of $m$. By Prop. 4.4 of \cite{SAS}, it follows that $\{1_{X+(m,n)}:(m,n) \in \bbz^2\}$ separates points of $Y$. Thanks to Thm. 3.3.5 of \cite{Arveson_invitation}, it follows that the translates of $X$ generate the $\sigma$-algebra of $Y$. 
Now $(2)$ follows.

For $k \geq 0$, let $E_k:=\{x \in X: \phi(x)+1=k\}$, and set $p_k=1_{E_k} \in L^{\infty}(\mathcal{D})$. Let $U$ be the bilateral shift on $\ell^2(\bbz)$. Let $S:L^{\infty}(\mathcal{D}) \to L^{\infty}(\mathcal{D})$ be defined by
\[Sf(x)=f(S^{-1}x).\]  

 Write  $L^{\infty}(Y)\rtimes \bbz^2$  as an iterated crossed product $(L^{\infty}(Y)\rtimes \bbz e_1)\rtimes \bbz e_2$. Note that $L^{\infty}(Y) \cong L^{\infty}(\mathcal{D})\otimes \ell^{\infty}(\bbz)$, and the 
action of $e_1$ is by translation by $1$ on the second factor. Thus, 
\[
L^{\infty}(Y) \rtimes \bbz e_1 \cong L^{\infty}(\mathcal{D})\otimes B(\ell^2(\bbz)).\]
Once $L^{\infty}(Y) \rtimes \bbz e_2$ is identified with $L^{\infty}(\mathcal{D})\otimes B(\ell^2(\bbz))$, the action of $e_2$ on the von Neumann algebra $L^{\infty}(\mathcal{D})\otimes B(\ell^2(\bbz))$ is given by 
$Ad(R)\circ (S \otimes 1)$, where the unitary $R$ is given by $R:=\sum_{k=0}^\infty p_k \otimes U^k$.  Thus, the automorphism corresponding to the $e_2$-action is outer conjugate to the automorphism $S \otimes 1$. 
Therefore, 
\[
L^{\infty}(Y)\rtimes \bbz^2 \cong (L^{\infty}(\mathcal{D})\rtimes_S \bbz) \otimes B(\ell^2(\bbz)).\]
It is well known that $L^{\infty}(\mathcal{D})\rtimes_S \bbz$ is a type $III_\lambda$ factor. Thus, $L^{\infty}(Y)\rtimes \bbz^2$ is a type $III_\lambda$ factor. The proof is complete. \hfill $\Box$

\begin{ppsn}
\label{type III examples}
Let $\beta>0$, $n$ be a positive integer, and let $\lambda:=e^{-\beta n}$. Then, there exists a pure $(\bbz^2,\bbn^2)$-space $(Y,X)$ and an $e^{-\beta c}$-conformal measure $\mu$ on $Y$ such that 
\begin{enumerate}
\item[(1)] the translates of $1_X$ generate $L^{\infty}(Y,\mu)$,
\item[(2)] $0<\mu(X)<\infty$, and
\item[(3)] the von Neumann algebra $L^{\infty}(Y,\mu) \rtimes \bbz^2$ is a factor of type $III_\lambda$.
\end{enumerate}
\end{ppsn}
\textit{Proof.} We apply Prop. \ref{Base type III} for the value $n\beta$ to get hold of a space $(Y_0,X_0)$ with a conformal measure. Thus, let $(Y_0,X_0)$ be a a pure $(\bbz^2,\bbn^2)$-space, and let $\nu$ be an $e^{-\beta nc}$-conformal measure on $Y_0$ such that 
\begin{enumerate}
 \item[(1)] the translates of $1_{X_0}$ generate $L^{\infty}(Y_0,\nu)$,
\item[(2)] $0<\nu(X_0)<\infty$, and
\item[(3)] the von Neumann algebra $L^{\infty}(Y_0)\rtimes \bbz^2$ is a factor of type $III_\lambda$.
\end{enumerate}
Denote the subgroup $n\bbz^2$ of $\bbz^2$ by $H$. As $H$ is isomorphic to $\bbz^2$, we can treat $Y_0$ as a $H$-space. We induce the action of $H$ on $Y_0$ to get a $\bbz^2$-space $Y$. 

Let $I:=\{0,1,2,\cdots,n-1\}$. 
Let $Y:=Y_0 \times I \times I$. Define a $\bbz^2$-action on $Y$ as follows:
\begin{equation*}
(y,p,q)+e_1:=\begin{cases}
 (y,p+1,q)  & \mbox{ if
} p < n-1,\cr
   &\cr
   (y+e_1,0,q) &  \mbox{ if } p=n-1,
         \end{cases}
\end{equation*}
and
\begin{equation*}
(y,p,q)+e_2:=\begin{cases}
 (y,p,q+1)  & \mbox{ if
} q < n-1,\cr
   &\cr
   (y+e_2,p,0) &  \mbox{ if } q=n-1.
         \end{cases}
\end{equation*}
Set $X:=X_0 \times I \times I$. It is clear that $X+\bbn^2 \subset X$. 
Note that for $r,s \geq 0$, we have
$X+(nr,ns)=(X_0+(r,s))\times I \times I$. 
Thus, 
\[
\bigcap_{(r,s) \in \bbn^2}(X+(nr,ns))=\Big(\bigcap_{(r,s) \in \bbn^2}X_0+(r,s)\Big)\times I \times I=\emptyset.\]
Also, for $r,s \geq 0$, 
$X-(nr,ns)=(X_0-(r,s))\times I \times I$. 
Consequently,
\[
\bigcup_{(r,s) \in \bbn^2}(X-(nr,ns))=\Big (\bigcup_{(r,s) \in \bbn^2}X_0-(r,s)\Big)\times I \times I=Y_0\times I \times I.\]
Thus, $(Y,X)$ is a pure $(\bbz^2,\bbn^2)$-space. Define a measure $\mu$ on $Y$ by setting
\[
d\mu:=e^{-\beta p}e^{-\beta q}d\nu dp dq.\]
Here, $dp dq$ is the counting measure on $I \times I$. Using the fact that $\nu$ is $e^{-\beta nc}$-conformal, it is routine to check that $\mu$ is $e^{-\beta c}$-conformal. We leave this verification
to the reader. Clearly, $\mu(X) \in (0,\infty)$. 

Let $U:=\{U_{(r,s)}\}_{(r,s) \in \bbz^2}$ be the Koopman representation of $\bbz^2$ on $L^{2}(Y_0,\nu)$. Denote the Hilbert space $L^{2}(X_0,\nu)$ by $H$, and let $V:=\{V_{(r,s)}\}_{(r,s) \in \bbn^2}$ be the isometric representation of $\bbn^2$ on $H$ defined by 
\[
V_{(r,s)}:=U_{(r,s)}|_{H}.\]

Let $\widetilde{U}:=\{\widetilde{U}_{(r,s)}\}_{(r,s) \in \bbz^2}$ be the Koopman representation of $\bbz^2$ on $L^{2}(Y,\mu)$. Denote the Hilbert space $L^{2}(X,\mu)$ by $\widetilde{H}$, and let $\widetilde{V}:=\{\widetilde{V}_{(r,s)}\}_{(r,s) \in \bbn^2}$ be the isometric representation
on $\widetilde{H}$ defined by 
\[
\widetilde{V}_{(r,s)}:=\widetilde{U}_{(r,s)}|_{\widetilde{H}}.\]

Let $\{\delta_p:p=0,1,2,\cdots,n-1\}$ be the standard orthonormal basis for $\ell^{2}(I)$. Let $U:\ell^{2}(I) \to \ell^2(I)$ be the unitary operator defined by
\begin{equation*}
U(\delta_p):=\begin{cases}
 \delta_{p+1}  & \mbox{ if
} p < n-1,\cr
   &\cr
   \delta_0 &  \mbox{ if } p=n-1.
         \end{cases}
\end{equation*}
Let $P$ be the orthogonal projection onto $\bbc \delta_{n-1}$. Note that we can identify $\widetilde{H}$ with $H \otimes \ell^{2}(I)\otimes \ell^{2}(I)$. Then, 
\begin{align*}
\widetilde{V}_{(1,0)}&=V_{(1,0)}\otimes UP \otimes 1+1 \otimes U(1-P)\otimes 1\\
\widetilde{V}_{(0,1)}&=V_{(0,1)}\otimes 1 \otimes UP+1 \otimes 1 \otimes U(1-P).
\end{align*}

We claim that $\widetilde{V}$ is irreducible. Let $M$ be the von Neumann algebra generated by $\{\widetilde{V}_{(r,s)}:(r,s) \in \bbn^2\}$. Suppose $T \in M^{'}$. 
Note that for $(r,s) \in \bbn^2$, 
\[
\widetilde{V}_{(nr,ns)}=V_{(r,s)}\otimes 1 \otimes 1.\]
Since $V$ is irreducible, it follows that there exists $R \in B(\ell^2(I)\otimes \ell^2(I))$ such that  $T=1\otimes R$.

The fact that $T\widetilde{V}_{(1,0)}=\widetilde{V}_{(1,0)}T$ implies that 
\[
V_{(1,0)}\otimes R(UP \otimes 1)+1 \otimes R(U(1-P)\otimes 1)=V_{(1,0)}\otimes (UP \otimes 1)R+1 \otimes (U(1-P)\otimes 1)R.\]

Since $V_{(1,0)}$ is an isometry which is not a unitary (as it has an eigenvalue of modulus strictly less than one), we have 
\begin{align*}
R(UP \otimes 1)&=(UP \otimes 1)R\\
R(U(1-P)\otimes 1)&=(U(1-P)\otimes 1)R.
\end{align*}
Routine simplifications of the above equations imply that  $R$ commutes with $U \otimes 1$ and $P \otimes 1$. Clearly, $\{U\otimes 1,P \otimes 1\}$ generates $B(\ell^2(I))\otimes 1$. Thus, $R \in 1 \otimes B(\ell^2(I))$. 

Similarly, the fact that $T\widetilde{V}_{(0,1)}=\widetilde{V}_{(0,1)}T$ leads to the conclusion that $R$ commutes with $1\otimes U$ and $1\otimes P$. Hence, $R \in B(\ell^2(I))\otimes 1$. Consequently, $R$ is a scalar multiple of the identity.
Hence, $\widetilde{V}$ is irreducible. It follows from Prop. \ref{trick to produce type III lambda} that the translates of $1_X$ generate $L^{\infty}(Y,\mu)$. 

Denote the subgroup $n\bbz^2$ of $\bbz^2$ by $H$, and let $G:=\bbz^2$.  Identify $H$ with $\mathbb{Z}^{2}$ via the map $H \ni (r,s) \to \frac{1}{n}(r,s) \in \mathbb{Z}^{2}$. Once $H$ is identified with $\bbz^2$, we can view $Y_0$ as a $H$-space. 
The $G$-space $Y$ is clearly the one obtained by inducing the $H$-space $Y_0$, i.e.
\[
Y:=Y_0 \times_H  G.\]
It can be proved  that $L^{\infty}(Y_0 \times_H G) \rtimes G \cong (L^{\infty}(Y_0) \rtimes H)\otimes B(\ell^2(G/H))$. Hence, as  $L^{\infty}(Y_0) \rtimes H$ is a factor of type $III_\lambda$,  $L^{\infty}(Y)\rtimes G$ is a factor of type $III_\lambda$. \hfill $\Box$

Theorem \ref{main1} is now a consequence of Prop. \ref{type III examples}, Prop. \ref{trick to produce type III lambda}, Prop. \ref{reduction to 1}, Prop. \ref{symmetry} and Prop. \ref{bijection between conformal measures}.

\bibliography{references}
 \bibliographystyle{amsplain}

 \vspace{2.5 mm}
 
 \noindent
{\sc S. Sundar}
(\texttt{sundarsobers@gmail.com})\\
         {\footnotesize  Institute of Mathematical Sciences, \\
         A CI of Homi Bhabha National Institute, \\
CIT Campus, Taramani, Chennai, 600113, \\Tamilnadu, INDIA.}\\

 \end{document}